\documentclass{gtart}

\def\ifplaintex{\expandafter\ifx\csname documentclass\endcsname\relax}

\def\gtp{{\mathsurround=0pt\it $\cal G\mskip-2mu$eometry \&\ 
$\cal T\!\!$opology $\cal P\!$ublications}}  

\def\recd{{\small Received:\qua\receiveddate\ifx\reviseddate\relax
\else\qquad Revised:\qua\reviseddate\fi\par}} 

\def\lognumber#1{\def\thelognumber{#1}}
\def\volumenumber#1{\def\thevolumenumber{#1}}
\def\volumeyear#1{\def\thevolumeyear{#1}}
\def\papernumber#1{\def\thepapernumber{#1}}
\def\pagenumbers#1#2{\def\startpage{#1}\def\finishpage{#2}}
\def\published#1{\def\publishdate{#1}}

\def\received#1{\def\receiveddate{#1}}
\def\revised#1{\def\reviseddate{#1}}
\def\accepted#1{\def\accepteddate{#1}}

\def\asciiaddress#1{\def\theasciiaddress{#1}}
\def\asciiemail#1{\def\theasciiemail{#1}}

\long\def\asciiabstract#1{\long\def\theasciiabstract{#1}}

\let\\\par\let\thelognumber\relax\let\thevolumenumber\relax
\let\thepapernumber\relax\let\thevolumeyear\relax\let\startpage\relax
\let\finishpage\relax\let\publishdate\relax\let\receiveddate\relax
\let\reviseddate\relax\let\accepteddate\relax\let\theasciititle\relax
\let\theasciiauthors\relax\let\theasciiaddress\relax
\let\theasciiabstract\relax

\let\theasciiemail\relax

\ifplaintex
\font\logobig=cmssbx10 scaled 3836
\font\logomed=cmssbx10 scaled 2557
\else
\font\logobig=cmssbx10 scaled 4200
\font\logomed=cmssbx10 scaled 2800
\fi

\long\def\makeagttitle{   
\count0=\startpage
\agt\hfill      
\hbox to 45truept{\vbox to 0pt{\vglue -13truept{\logomed A\kern -.37em{\logobig 
T}\kern -.38em G}\vss}\hss}
\break
{\small Volume \thevolumenumber\ (\thevolumeyear)
\startpage--\finishpage\nl
Published: \publishdate}

\vglue .25truein

{\parskip=0pt\leftskip 0pt plus
1fil\def\\{\par\smallskip}{\Large\bf\thetitle}\par\medskip} \vglue
0.05truein

{\parskip=0pt\leftskip 0pt plus 1fil\def\\{\par}{\sc\theauthors}
\par\medskip}%
 
\vglue 0.03truein 

{\small\leftskip 25truept\rightskip 25truept{\bf Abstract}\stdspace\theabstract

{\bf AMS Classification}\stdspace\theprimaryclass
\ifx\thesecondaryclass\relax\else; \thesecondaryclass\fi\par
{\bf Keywords}\stdspace \thekeywords\par}\vglue 7truept

}   

\ifplaintex
\font\phead=cmsl9 scaled 950
\font\pnum=cmbx10 scaled 913
\font\pfoot=cmsl9 scaled 950
\headline{\vbox to 0pt{\vskip -4.5mm\line{\small\phead\ifnum
\count0=\startpage ISSN 1472-2739 (on-line) 1472-2747 (printed)
\hfill {\pnum\folio}\else\ifodd\count0\def\\{ }%
\ifx\theshorttitle\relax\thetitle\else\theshorttitle\fi\hfill{\pnum\folio}
\else\def\\{ and }{\pnum\folio}\hfill\ifx\theshortauthors\relax\theauthors
\else\theshortauthors\fi\fi\fi}\vss}}
\footline{\vbox to 0pt{\vglue 0mm\line{\small\pfoot\ifnum\count0=\startpage
\copyright\ \gtp\hfill\else
\agt, Volume \thevolumenumber\ (\thevolumeyear)\hfill\fi}\vss}}
\else
\font=cmsl9 scaled 1050
\font\lnum=cmbx10 
\font=cmsl9 scaled 1050
\makeatletter
\def\@oddhead{{\small\lhead\ifnum\count0=\startpage ISSN 1472-2739 
(on-line) 1472-2747 (printed)\hfill {\lnum\number\count0}\else\ifodd\count0
\def\\{ }\ifx\theshorttitle\relax \thetitle \else\theshorttitle\fi\hfill
{\lnum\number\count0}\else\def\\{ and }{\lnum\number\count0}
\hfill\ifx\theshortauthors\relax 
\theauthors\else\theshortauthors\fi\fi\fi}}\def\@evenhead{\@oddhead}
\def\@oddfoot{\small\lfoot\ifnum\count0=\startpage\copyright\ \gtp\hfill\else
\agt, Volume \thevolumenumber\ (\thevolumeyear)\hfill\fi}
\def\@evenfoot{\@oddfoot}
\makeatother
\fi
\let\maketitlepage\makeagttitle
\let\makeshorttitle\maketitlepage
\let\maketitle\maketitlepage

\newwrite\gtoutfile
\long\gdef\makeheadfile{  
{\def\\{, }\def\s{ }
\immediate\openout\gtoutfile head.xxx
\immediate\write\gtoutfile{Proxy-for: \ifx\theasciiauthors\relax
\theauthors\else\theasciiauthors\fi\s<\ifx\theasciiemail\relax\theemail\else\theasciiemail\fi>}
\immediate\write\gtoutfile{\noexpand\\}
\immediate\write\gtoutfile{Authors: \ifx\theasciiauthors\relax
\theauthors\else\theasciiauthors\fi}
{\def\\{ }\immediate\write\gtoutfile{Title: \ifx\theasciititle\relax
\thetitle\else\theasciititle\fi}}
\immediate\write\gtoutfile{Subj-class: GT or SG, GR etc}
\immediate\write\gtoutfile{MSC-class: \theprimaryclass\ifx\thesecondaryclass\relax\else, \thesecondaryclass\fi}
\immediate\write\gtoutfile{Journal-ref: Algebr. Geom. Topol. \thevolumenumber\s
(\thevolumeyear) \startpage-\finishpage}
\immediate\write\gtoutfile{Comments: Published by Algebraic and
Geometric Topology at}
\immediate\write\gtoutfile{\s\s\s  http://www.maths.warwick.ac.uk/agt/AGTVol\thevolumenumber/agt-\thevolumenumber-\thepapernumber.abs.html}
\immediate\write\gtoutfile{\noexpand\\}
\immediate\write\gtoutfile{}
\ifx\theasciiabstract\relax
\immediate\write\gtoutfile{\theabstract}\else
\immediate\write\gtoutfile{\theasciiabstract}\fi
\immediate\write\gtoutfile{}
\immediate\write\gtoutfile{\noexpand\\}
\immediate\write\gtoutfile{}
\immediate\closeout\gtoutfile}}  

\def\maketitlepage{\makeagttitle\makeheadfile}
\let\makeshorttitle\maketitlepage
\let\maketitle\maketitlepage

\lognumber{42}
\volumenumber{3}
\volumeyear{2003}
\papernumber{42}
\published{12 December 2003}
\pagenumbers{1167}{1224}
\received{5 February 2003}
\revised{1 December 2003}
\accepted{11 December 2003}

\usepackage{amssymb, amsmath, latexsym, amscd, graphicx, subfigure}

\newcommand{\lbl}{\label}

\newtheorem{theorem}{Theorem}
\newtheorem{corollary}{Corollary}
\newtheorem{lemma}{Lemma}
\newtheorem{proposition}{Proposition}
\theoremstyle{definition}

\newtheorem{definition}{Definition}

\newtheorem{remark}{Remark}

\def\bdry{\partial}
\def\bE{\bdry_E}

\newcommand{\F}{\mathcal F}
\newcommand{\G}{X}

\newcommand{\R}{\mathbb R}

\newcommand{\g}{\mathcal G}

\newcommand{\sy}{\mathfrak {sp}}
\renewcommand{\a}{\mathcal L}
\newcommand{\OX}{U}             
\newcommand{\OY}{IU}            
\newcommand{\CX}{\widehat{\OX}} 
\newcommand{\CY}{\widehat{\OY}} 
\newcommand{\s}{A}              
\newcommand{\V}{\widehat{f\mathcal  G}}      
\newcommand{\C}{f\mathcal G}               
\newcommand{\W}{\widehat{fr\mathcal G}}    
\newcommand{\fr}{fr\mathcal G}              

\def\O{ {\mathcal O} }          
\def\Sp{ {\mathcal {OS}} }        
\def\G{ {\mathcal {OG}} }         
\def\A{ {\mathcal {LO}}_n }       
\def\R{{\mathbb R}}             
\def\ss{\mathbb S}                
\def\is{\Lambda}              
\def\X{{\bf X}}     

\begin{document}

\title{On a theorem of Kontsevich}
\authors{James Conant\\Karen Vogtmann}
\address{Department of Mathematics, University of Tennessee\\
         Knoxville, TN 37996, USA}
\email{jconant@math.utk.edu}
\secondaddress{Department of Mathematics, Cornell University\\
         Ithaca, NY 14853-4201, USA}
\secondemail{vogtmann@math.cornell.edu}
\asciiemail{jconant@math.utk.edu, vogtmann@math.cornell.edu}

\asciiaddress{Department of Mathematics, University of 
Tennessee\\Knoxville, TN 37996, USA\\and\\Department of 
Mathematics, Cornell University\\Ithaca, NY 14853-4201, USA}

\begin{abstract}
In \cite{K1} and \cite{K2}, M. Kontsevich introduced graph
homology as a tool to compute the homology of three infinite dimensional
Lie algebras, associated to the three operads ``commutative,"
``associative" and ``Lie." We generalize his theorem to all cyclic
operads, in the process giving a more careful treatment of the
construction than in Kontsevich's original papers. We also give a more
explicit treatment of the isomorphisms of graph homologies with the
homology of moduli space and $Out(F_r)$   outlined by Kontsevich.  In
\cite{CoV} we defined a Lie bracket and cobracket on the commutative
graph complex, which was extended in \cite{C} to the case of all cyclic
operads.  These operations form a   Lie bi-algebra on a
natural subcomplex.   We show that in the associative and Lie cases the
subcomplex on which the bi-algebra structure exists carries all of the
homology, and we explain why the subcomplex in the commutative case does
not.   
\end{abstract}

\asciiabstract{In [`Formal (non)commutative symplectic geometry', The
Gelfand Mathematical Seminars (1990-1992) 173-187, and `Feynman
diagrams and low-dimensional topology', First European Congress of
Mathematics, Vol. II Paris (1992) 97--121] M. Kontsevich introduced
graph homology as a tool to compute the homology of three infinite
dimensional Lie algebras, associated to the three operads
`commutative,' `associative' and `Lie.' We generalize his theorem to
all cyclic operads, in the process giving a more careful treatment of
the construction than in Kontsevich's original papers. We also give a
more explicit treatment of the isomorphisms of graph homologies with
the homology of moduli space and Out(F_r) outlined by Kontsevich.  In
[`Infinitesimal operations on chain complexes of graphs',
Mathematische Annalen, 327 (2003) 545-573] we defined a Lie bracket
and cobracket on the commutative graph complex, which was extended in
[James Conant, `Fusion and fission in graph complexes', Pac. J. 209
(2003), 219-230] to the case of all cyclic operads.  These operations
form a Lie bi-algebra on a natural subcomplex.  We show that in the
associative and Lie cases the subcomplex on which the bi-algebra
structure exists carries all of the homology, and we explain why the
subcomplex in the commutative case does not.}

\maketitle

\section{Introduction}

In the papers \cite{K1} and \cite{K2} M. Kontsevich sketched an elegant theory which  relates the homology of certain infinite-dimensional
 Lie algebras to various invariants in low-dimensional topology and group theory.  The infinite-dimensional Lie
algebras arise as generalizations of the Lie algebra of polynomial functions on $\mathbb R^{2n}$  under the classical Poisson bracket or, equivalently,
the Lie algebra of polynomial vector fields on $\mathbb R^{2n}$ under the Lie bracket.  One thinks of $\mathbb R^{2n}$ as a symplectic manifold, and notes
that these Lie algebras each contain a copy of the symplectic Lie algebra $\sy(2n)$. The relation with  topology and group theory is
established by interpreting the $\sy(2n)$-invariants in the exterior algebra of the Lie algebra in terms of graphs, and then exploiting both new
and established connections between graphs and areas of low-dimensional topology and group theory.  These connections include the construction of
3-manifold and knot invariants via data associated to trivalent graphs, the study of automorphism groups of free groups using the space of marked metric
graphs (Outer space), and the use of ribbon graph spaces to study mapping class groups of punctured surfaces.

This paper is the outcome of a seminar held at Cornell, organized by the second author, devoted to understanding Kontsevich's theory.  
Kontsevich describes three variations of his theory, in the commutative, associative and Lie ``worlds." 
Kontsevich's papers skip many definitions and details, and, as we discovered, have
a gap in the proof of the main theorem relating symplectic invariants and graph homology. In this paper we explain Kontsevich's theorem carefully, in the
more general setting of cyclic operads. In particular, we adapt a fix that Kontsevich communicated to us for the commutative case to the general case. 
We then specialize to the  Lie, associative and commutative operads, which are the three worlds which Kontsevich considered in his original
papers.  Using a filtration of Outer space
indicated by Kontsevich, we show that the primitive part of the homology of the Lie graph  complex is the direct sum of the cohomologies of
$Out(F_r)$, and the primitive part of the homology of the associative graph complex is the direct sum of the cohomologies
of moduli spaces   (or equivalently mapping class groups) of punctured surfaces.   We then recall the
Lie bracket and cobracket which
 we defined on the commutative
graph complex in \cite{CoV}, and which was extended in \cite{C} to the case of all cyclic
operads.  These operations form a bi-algebra structure on a natural
subcomplex.   We show that in the associative and Lie cases the
subcomplex on which  the bi-algebra structure exists carries all of the
homology, and we explain why the subcomplex in the commutative case does
not.

  Many people   contributed to 
this project. Participants in the Cornell seminar included David Brown, Dan Ciobutaru,
Ferenc Gerlits, Matt Horak, Swapneel Mahajan and  Fernando Schwarz.  We are particularly indebted to
Ferenc Gerlits and Swapneel Mahajan, who continued to work with us on understanding these papers after
the scheduled seminar was over.  Mahajan has written a separate exposition of some of the material here,
with more information about the relation with classical symplectic geometry, using  what he calls
\emph{reversible} operads, which are closely related to cyclic operads. His treatment of the theorem
relating graphs and Lie algebra invariants uses Kontsevich's original fix, and therefore does not work
for every cyclic operad. In particular, it cannot handle the Lie case. 
 A succinct outline of Kontsevich's theory can  be found in the thesis of Gerlits, which also includes a
careful study of the Euler characteristic using Feynman integrals.   We hope that these different expositions, with different emphases, will help to
make Kontsevich's beautiful theory more broadly accessible.  

The present paper is organized as follows.  In Section 2 we develop  the theory for general cyclic operads.  After reviewing the definition
of cyclic operad we define a chain complex parameterized by graphs whose vertices are   ``colored" by
operad elements.  This chain complex was introduced in a more general setting by Getzler and Kapranov
in \cite{GeKa2}, and was studied by Markl in 
\cite{Markl}.   We then construct a functor from 
 cyclic operads to symplectic Lie algebras, as the direct limit of functors indexed by natural numbers.
 We then show how to use invariant theory of the symplectic Lie algebra to 
define a map from the Chevalley-Eilenberg complex of the Lie algebra to the above chain complex of graphs which is an isomorphism on
homology. 
   
In section 3, we specialize to the Lie operad.  By using a
filtration of Outer space indicated by Kontsevich, we prove that the primitive (connected) part of the graph complex computes the cohomology of the
groups
$Out(F_r)$ of outer automorphisms of finitely-generated free groups. We also  prove that inclusion of the subcomplex spanned by 1-particle
irreducible graphs (i.e.\ graphs with no separating edges) is an isomorphism on homology. This is of interest because, as was shown in \cite{CoV, C} this subcomplex carries a graded Lie bi-algebra
structure; this implies, among other things, that the homology of
  the groups $Out(F_r)$ is the primitive part of a differential graded Hopf algebra. 
In section 4 we note how the theory applies in the case of the associative operad.   In this case, the primitive part of the graph complex is shown to
compute the cohomology of mapping class groups of punctured surfaces.  The proof proceeds by restricting the filtration of Outer space to the 
``ribbon graph" subcomplexes, on which  mapping class groups of  punctured surfaces act.   As in the Lie case, we show that the subcomplex of 1-particle
irreducible graphs carries all of the homology, so that the direct sum of homologies of mapping class groups is the primitive part of a differential
graded Hopf algebra.
Finally, in section 5 we reconsider the commutative case, which was the focus of our paper \cite{CoV}.  We give a geometric description of 
the primitive part of commutative graph homology, as  the relative homology of a completion of Outer space modulo the subcomplex at
infinity, with certain twisted coefficients.   This relative homology measures the difference between the relative quotient of Outer space by $Out(F_r)$ and the quotient by
the subgroup
$SOut(F_r)$ of outer automorphisms which map to $SL(r,\mathbb Z)$ under the natural map $Out(F_r)\to GL(r,\mathbb Z)$.  Using this geometric description of
graph homology, we explain why the one-particle irreducible subcomplex of the graph complex does not have the same homology as the full
complex, unlike the Lie and associative cases.

Acknowledgments: In addition to the seminar participants mentioned above, we would like to thank Martin Markl, Steve Shnider and Jim Stasheff for pointing out several typos and unclear statements. Finally, the first author was partially supported by NSF VIGRE grant  DMS-9983660, and the second author was partially supported by NSF grant DMS-0204185.

\section{The general case - Cyclic Operads}

\subsection{Review of cyclic operads} Throughout the paper we will work in the category of real vector spaces.  
Thus an \emph{operad} $\O$ is a collection of real
vector spaces $\O[m]$, $m\geq 1$, with 
\begin{itemize}
\item[1)] a composition law:
$$\gamma\colon\O[m]\otimes \O[i_1]\otimes \ldots\otimes \O[i_m]\to  
\O[i_1+\ldots+i_m],$$
\item[2)] a right action of the symmetric group $\Sigma_m$ on $\O[m]$, and
\item[3)] a unit $1_\O\in \O[1]$
\end{itemize}
  satisfying appropriate axioms governing the unit, associativity, and
$\Sigma_m$-equiv\-ar\-iance of the composition law (see \cite{MSS} for a
complete list of axioms).  Intuitively, an element of $\O[m]$ is an
object with $m$ numbered input slots and an ouput slot. The symmetric
group acts by permuting the numbering of the input slots. The
composition $\gamma(o\otimes o_1\otimes\ldots\otimes o_m)$ plugs the
output of $o_j\in \O[i_j]$ into the $j$th input slot of $o$,
renumbering the input slots of the result consistently (see
Figure~\ref{composition}).
\begin{figure}[ht!]
\begin{center}
\includegraphics[width=.9\hsize]{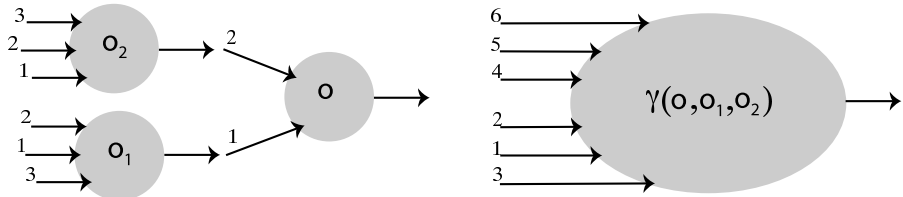}
\caption{Operad composition}\lbl{composition}
\end{center}
\end{figure}
If $o_i=1_{\O}$  for $i\neq k$ and $o_k=o'$, we call the result the $k$th composition of $o'$ with $o$.

A \emph{cyclic operad} is an operad where the action of the symmetric group $\Sigma_m$
extends to an action of $\Sigma_{m+1}$ in a way compatible with the axioms. This concept was introduced in
\cite{GeKa} (see also \cite{MSS}).
The intuitive idea  is that in a cyclic operad, the output slot can also serve as an input slot, and any input slot
can serve as the output slot.
Thus, if we number the output slot as well (say with $0$), then
$\Sigma_{m+1}$ acts by permuting the numbers on all input/output slots.  Modulo the $\Sigma_{m+1}$ action we can
compose two elements using one input/output slot of each. 

\begin{figure}[ht!]
\begin{center}
\subfigure[Labeled $6$-star]{\includegraphics[width=.27\linewidth]{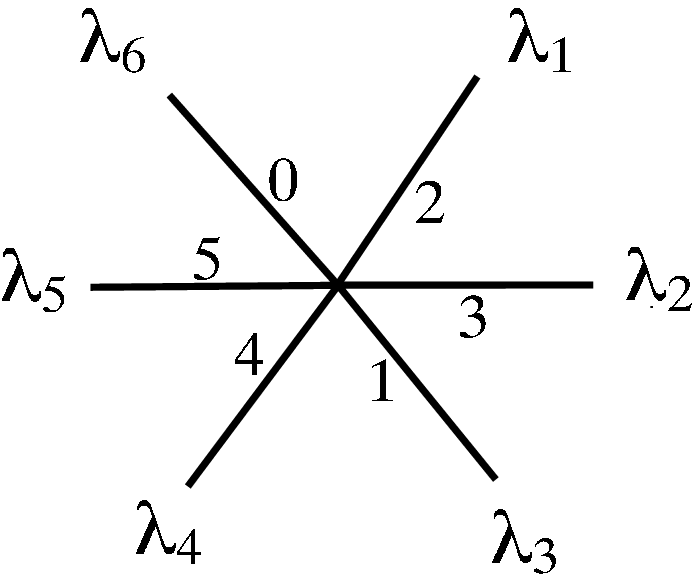}}
\space\space\space\space\space
\subfigure[Superimposing an operad element and an $m$-star]{\includegraphics[width=.3\linewidth]{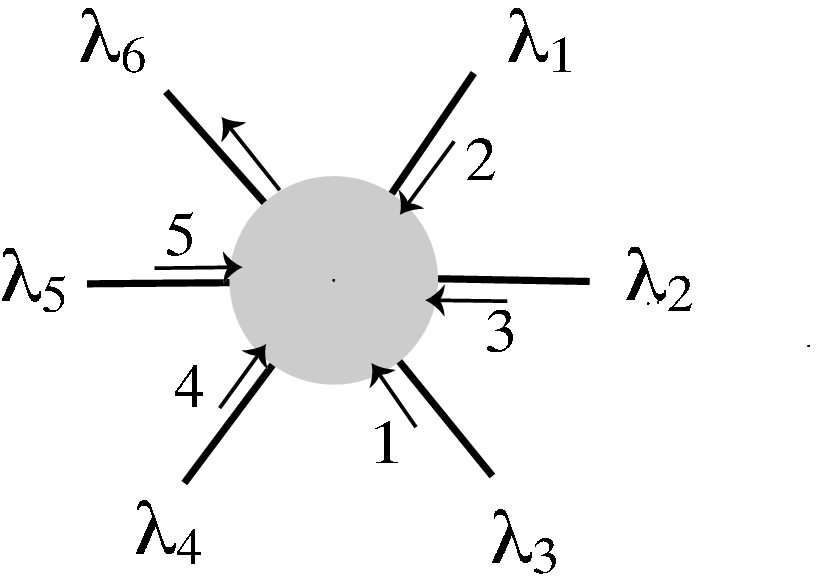}}
\space\space\space\space\space
\subfigure[$\O$-spider]{\includegraphics[width=.27\linewidth]{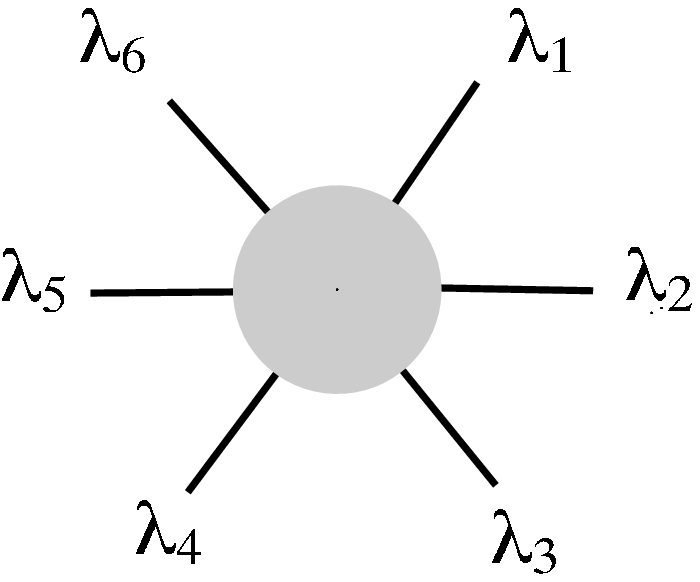}}
\caption{``Coloring" an $m$-star with an element of $\O$ to make an $\O$-spider}\lbl{spider}
\end{center}
\end{figure}

We can bring the actual situation closer to the intuition as follows.  For each integer $m\geq 2$, let 
$*_m$ be the
\emph{$m$-star}, the unrooted tree with one internal vertex and $m$ leaves $\lambda_1,\ldots,\lambda_m$. A \emph{labeling} of $*_m$ is, by
definition, a bijection from  the set of leaves to the numbers 
$0,\ldots, m-1$. We represent a labeling $L$ by placing the number $L(\lambda_i)$ on the leaf $\lambda_i$, 
close to the internal vertex of $*_m$. Figure~\ref{spider}a shows a labeled $6$-star.

  The symmetric group $\Sigma_{m}$ acts on the set of labelings of $*_m$, and we make the following
definition.

\begin{definition} Let $\O$ be a cyclic operad,  $m\geq 2$ an integer, and $\mathcal L$  the set of labelings of the $m$-star.  The space  $\Sp[m]$  of
\emph{$\O$-spiders with $m$ legs} is defined to be the space of coinvariants
$$\Sp[m]=(\bigoplus_{\mathcal L}\O[m-1])_{\Sigma_{m}},$$ 
where $(\cdot)_L : \O[m-1]\to \bigoplus_{L}\O[m-1]$ is
the natural inclusion into the $L$ summand, and $\Sigma_m$ acts by $\sigma\cdot (o_L)=(\sigma\cdot o)_{\sigma\cdot L}$.
\end{definition}

Every element of $\Sp[m]$ is an equivalence class 
$[o_L]$ for some $o\in \O[m-1]$ and labeling $L$ of $*_m$.  To see this, note that 
$[(o_1)_{L_1}+(o_2)_{L_2}]=[(o_1)_{L_1}]+[(o_2)_{L_2}]=[(o_1)_{L_1}]+[(\sigma\cdot o_2)_{L_1}]=[(o_1+\sigma\cdot o_2)_{L_1}]$ for $\sigma\in
\Sigma_m$ with $\sigma\cdot L_2=L_1$. We can think of
$o$ as sitting on top of
$*_m$ so that the labeling of $*_m$ corresponds to the numbering of the input/output slots (Figure~\ref{spider}b).

Modding out by the action of $\Sigma_{m}$ erases the labeling and the distinction between input and output slots
(Figure~\ref{spider}c). 
The picture explains the arachnoid terminology for elements of $\Sp[m]$.
\begin{figure}[ht!]
\begin{center}
\subfigure[Spiders to be mated using legs $\lambda$ and $\mu$]{\includegraphics[height=2.8cm]{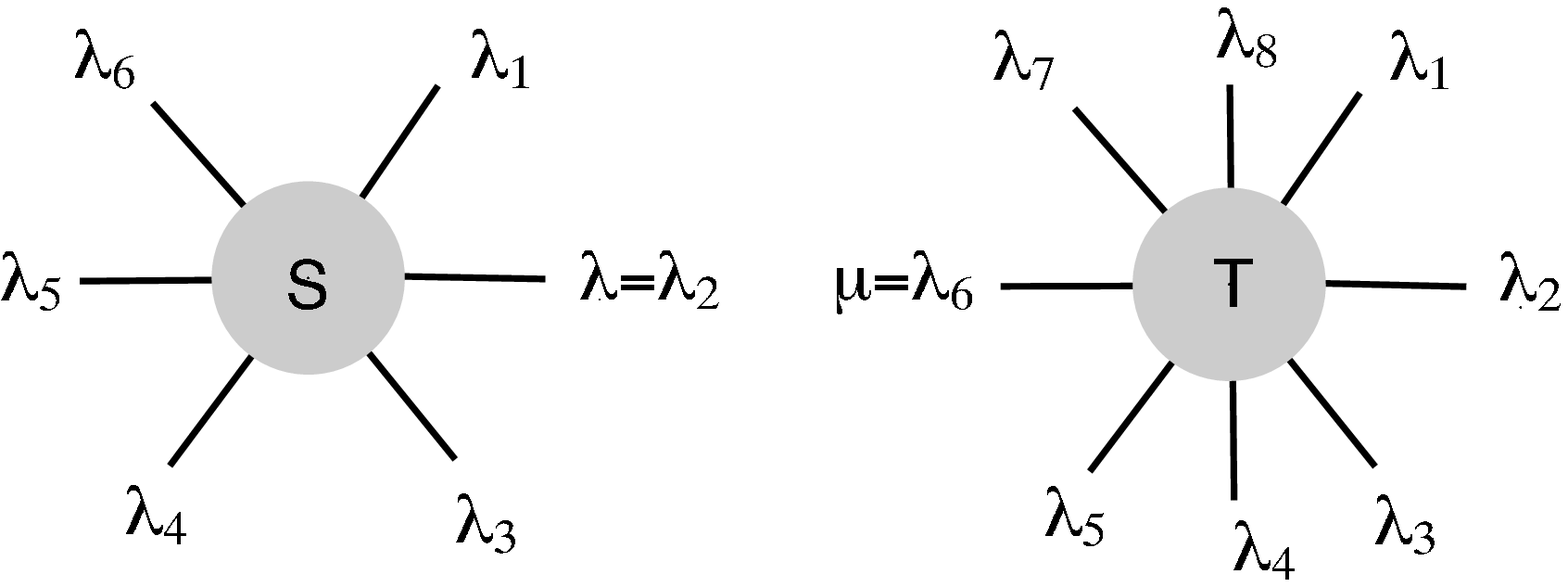}}
\hfill
\subfigure[Phase one]{\includegraphics[height=3cm]{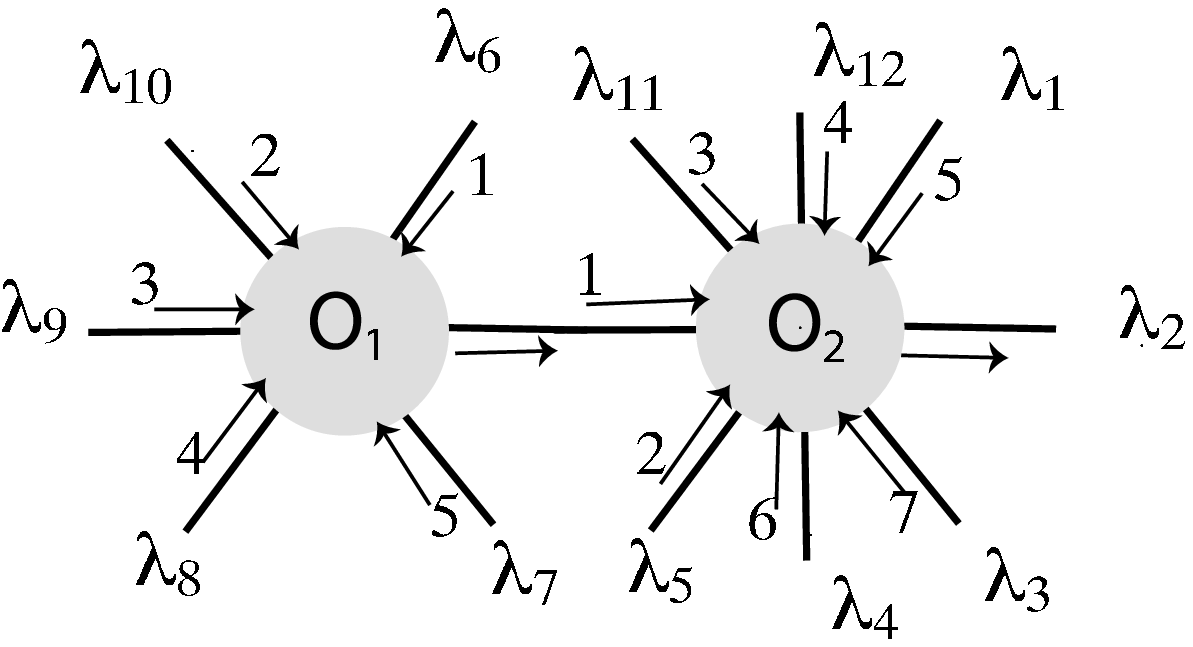}}
\subfigure[Phase two]{\includegraphics[height=3cm]{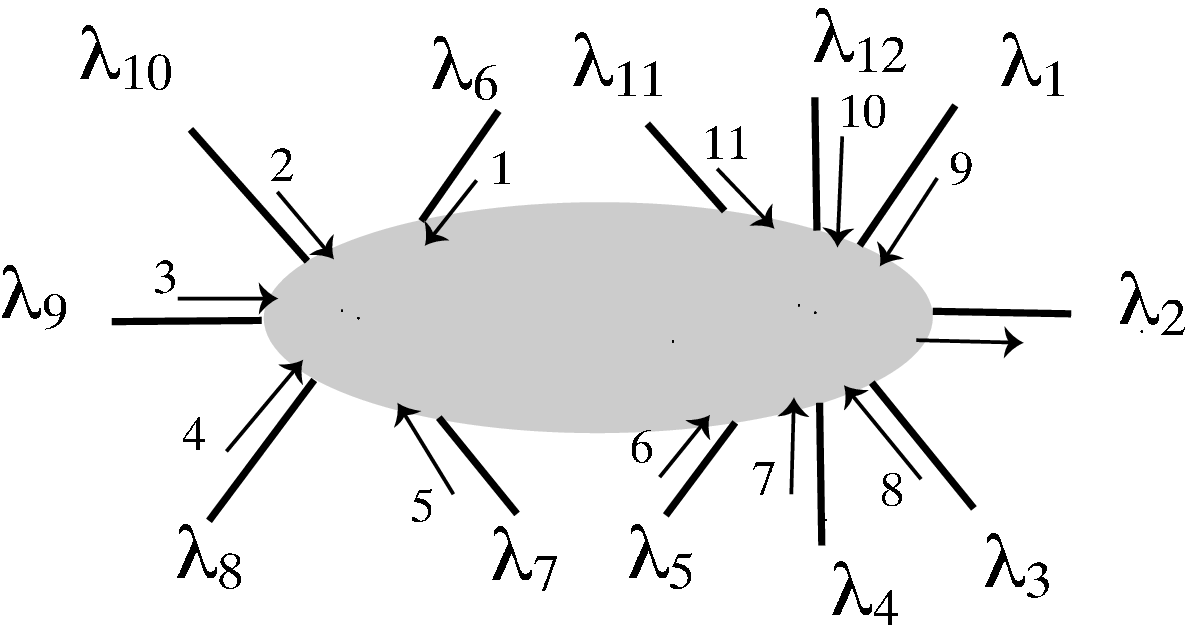}}
\hfill
\subfigure[Mated spiders]{\includegraphics[height=3cm]{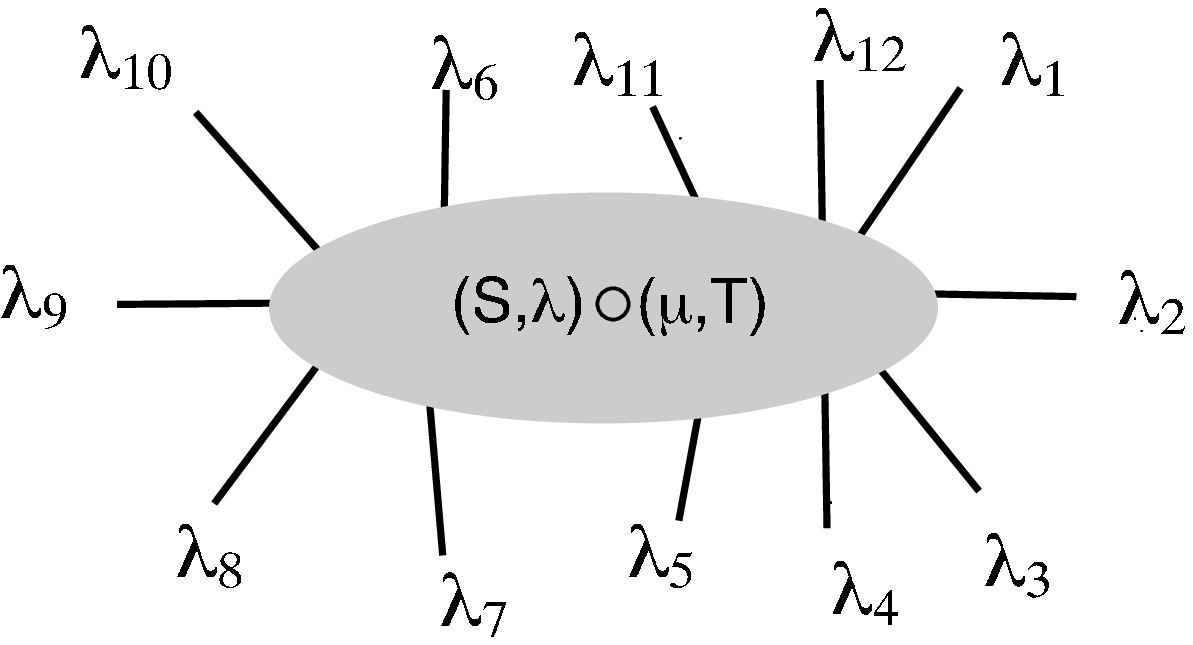}}
\caption{Mating spiders}\lbl{composed}
\end{center}
\end{figure}

The composition law
in $\O$ transforms to a \emph{mating law} in $\Sp=\bigoplus_{m\geq 2} \Sp[m]$, as follows.
Consider two $\O$-spiders, $S$ and $T$, and any pair of legs 
 $\lambda$ of $S$ and $\mu$ of $T$ (Figure~\ref{composed}a).

Choose representatives 
$(o_1)_{L_1}$ for $S$ and  $(o_2)_{L_2}$ for $T$ 
 such that $L_1(\lambda)=0$  (so that $\lambda$ corresponds to the output slot of $o_1$) and $L_2(\mu)=1$ (so $\mu$ corresponds to the first 
input slot of $o_2$).  
Connect the $\lambda,\mu$ legs together.
Rename the spider legs (other than $\lambda$ and $\mu$) so that the remaining legs of $S$ are inserted, in order, into the 
ordered set of legs of $T$, at the slot formerly occupied by  $\mu$. (Figure
~\ref{composed}b).   

Now contract the edge formed by $\lambda$ and $\mu$ to get an underlying $m$-star, and
compose the two operad elements along the corresponding input/output slots to obtain
$\gamma(o_2\otimes o_1 \otimes 1_\O\otimes\cdots\otimes 1_\O)_L$, where $L$ is the induced 
labeling (Figure~\ref{composed}c). 

The resulting equivalence class under the symmetric group action will be denoted by $(S,\lambda)\circ
(\mu,T)$ (Figure~\ref{composed}d).

\subsection{Examples}
We describe the three operads we will be focusing on: the commutative, Lie and associative operads. 
Figure~\ref{3spiders} shows examples of spiders in these operads. There are many other cyclic operads, 
for example the endomorphism operad and the Poisson operad. See \cite{GeKa} or \cite{MSS}. It is also worthwhile to note at this point that cyclic operads form a category, and there are obvious morphisms from the Lie operad to the associative operad, and from the associative to the commutative operad. See \cite{MSS}.
\begin{figure}[ht!] 
\begin{center}
\includegraphics[width=\linewidth]{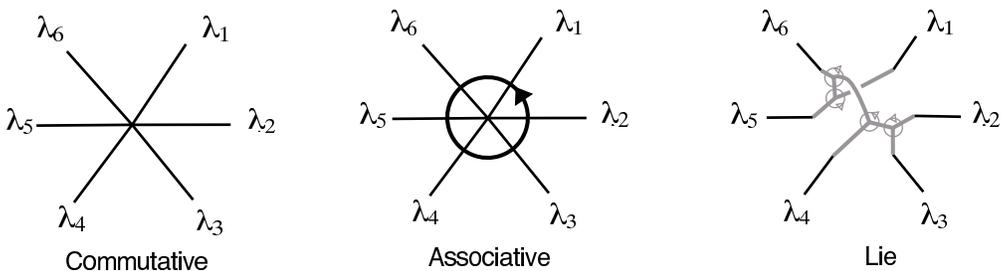}
\caption{Three types of spiders. In the Lie case, the picture is modulo IHX and AS
relations.}\lbl{3spiders} 
\end{center}
\end{figure}

\subsubsection{The commutative operad} 
  In the 
commutative operad, each $\mathcal O[m]$ is $1$-dimensional, with  trivial $\Sigma_m$ action.
The composition law is given by the canonical isomorphism (i.e.\ multiplication) $\mathbb R^{\otimes k}\cong \mathbb R$.
An $\O$-spider in this case is a copy of $*_m$, weighted by a real number. Mating is done by joining
legs $\lambda$ and
$\mu$ of two stars to form an edge, then contracting that edge to form a new star and multiplying the
weights.  

\subsubsection{The associative operad}
  In the associative operad,  each $\mathcal O[m]$ is spanned by rooted planar trees with one internal
vertex and $m$ numbered leaves. The   planar embedding that  each such tree comes with is equivalent to
a prescribed left-to-right ordering of the leaves. The symmetric group acts by permuting the numbers of
the leaves.  To compose two trees, we attach the root of the first tree to a leaf of the second tree,
then collapse the  internal edge we just created.  We order the leaves of the result so that the
leaves of the first tree are inserted, in order, at the position of the  chosen leaf in the second
tree. A basis element of the space $\Sp$ of $\O$-spiders is a copy of
$*_m$ with a fixed \emph{cyclic} ordering of the legs.  To mate two basic $\O$-spiders using legs $\lambda$ and $\mu$, we join $\lambda$ and $\mu$
to form a connected graph with one internal edge, then contract the internal edge. The cyclic order on the legs of each spider induces a cyclic
ordering of the legs of the spider which results from mating. Mating is extended linearly to all spiders, i.e.\ to spiders whose ``body" consists of a
linear combination of cyclic orderings.

\subsubsection{The Lie operad}  
  In the Lie operad, $\O[m]$ is the vector
space spanned by all
 rooted planar binary trees with $m$ numbered
leaves, modulo the subspace
spanned by all anti-symmetry and IHX relators.
Specifying a planar embedding of a tree is equivalent to giving a cyclic ordering of the edges adjacent to each interior vertex.
The anti-symmetry relation AS
says that switching the cyclic order at any vertex reverses the sign. The IHX (Jacobi) relation is well-known (see Figure
~\ref{IHX}).
\begin{figure}[ht!]
\begin{center}
\includegraphics[height=3cm]{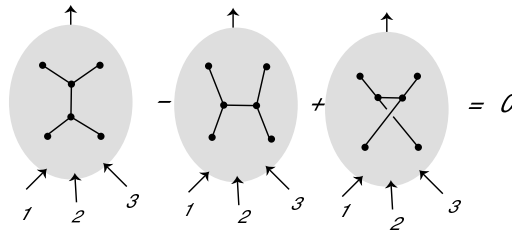}
\end{center}
\caption{The IHX relator}\lbl{IHX}
\end{figure}

The symmetric group $\Sigma_m$ acts by permuting the numbering of the leaves.   The composition rule in the operad attaches the root of
one tree to a leaf of another to form an interior edge of a new planar tree, then suitably renumbers
the remaining leaves. 

The space of $\O$-spiders, $\Sp$, is spanned by planar binary trees with $m$ numbered leaves, modulo AS and IHX,
but with no particular leaf designated as the root. Mating is accomplished by gluing two such trees together
at a leaf to form a single planar binary tree, then renumbering the remaining leaves suitably. 
Note that mating does not involve an edge collapse, as it did in the commutative and
associative cases.

\subsection{Graph homology of a cyclic operad $\O$}
In this section we construct a functor 
$$\text{\{Cyclic operads\}$\to$\{Chain complexes\}},$$ where the chain
complexes are spanned by oriented graphs with an element of $\Sp$ attached to each vertex. 
We
begin with a subsection discussing the appropriate notion of orientation on
graphs.  This subsection includes results which will be needed later when working with
specific operads; the reader interested
only in the basic construction can stop at the definition of orientation on first reading.

\subsubsection{Oriented graphs}\lbl{sec:orientation}

  By a {\it  graph} we mean a finite 1-dimensional CW complex. The set of
edges of a graph $X$ is denoted $E(X)$, the set of vertices $V(X)$ and the set of
half-edges $H(X)$. Let $H(e)$ denote the set of (two) half-edges contained in an edge $e$. 
There is an involution
$x\mapsto
\bar x$ on
$H(X)$, swapping the elements of $H(e)$ for each $e\in E(X)$. 

For an  $n$-dimensional vector space $V$, set 
 $\det(V):=\wedge^n V$.  An orientation on $V$ can be thought of as a unit vector in $det(V)$.
For a set $Z$, we denote by $\R Z$ the real vector space with basis Z.

\begin{definition} An \emph{ orientation} on
a graph
$X$   is a unit vector in
$$\det \R V(X) \otimes\bigotimes_{e\in E(X)}\det\R H(e).$$
\end{definition}

In other words, \emph{an orientation  on
$X$ is determined by ordering the vertices of $X$ and orienting each edge of
$X$}.  Two orientations are the same if they are obtained from one another by
an even number of edge-orientation switches and vertex-label swaps.

Our definition is different from the definition given in Kontsevich's papers
\cite{K1},\cite{K2} but, as we show below, it is equivalent for connected graphs.
(Note that Kontsevich defines his graph homology using only connected graphs.) 
  We follow ideas of Dylan Thurston
\cite{Thurston}, and begin
by recording a basic observation:

\begin{lemma}\lbl{exact} Let $0\to A\to B\to C\to D\to 0$ be an exact sequence of finite-dimensional
vector spaces.  Then there is a canonical isomorphism $$\det(A)\otimes \det(C)\to
\det(B)\otimes
\det(D).$$
\end{lemma}

 \begin{proof} 
 
For any short exact sequence $0\to U \stackrel{f}\to V\to W\to 0$   of
finite-dimensional vector spaces,  with 
$s\colon W\to V$   a splitting, the map 
$$\det(U)\otimes\det(W)\to \det(V)$$
 given by $u\otimes w\mapsto f(u)\wedge s(w)$ is an isomorphism, and is independent of the
choice of $s$.
The lemma now follows by splitting $0\to A\to B\to C\to D\to 0$ into two short exact
sequences.
\end{proof}

Kontsevich defines an orientation on a   graph to be an orientation of the vector
space $H_1(X;\mathbb R)\oplus \R{E(X)}$. The following proposition shows that this is equivalent
to our definition for connected graphs.

\begin{proposition}\lbl{orient} Let $X$ be a connected graph.  Then there is a canonical
isomorphism
 $$ \det(\R V(X))\otimes \bigotimes_e \det(\R
H(e)) \cong\det(H_1(X;\mathbb R))\otimes \det(\R E(X)).  $$
\end{proposition}

\proof 
For any graph $X$, we have an exact sequence
$0\to H_1(X;\mathbb R)\to C_1(X)$\break $\to C_0(E)\to H_0(X;\mathbb R)\to 0,$ so Lemma~\ref{exact} gives
  a canonical isomorphism $$\det(H_1(X;\mathbb R))\otimes \det(C_0(X))\cong \det(C_1(X))\otimes
\det(H_0(X;\mathbb R)).\eqno{(1)}$$

 $C_0(X)$ has a canonical basis consisting of the vertices of $X$, so that $C_0(X)$ can be
identified with
$\R V(X)$:
$$\det(C_0(X))\cong \det(\R V(X)).
$$In order to give a chain in  $C_1(X)$, on the other hand, you need to
prescribe orientations on all of the edges, so that
 $$\det(C_1(X))\cong \det(\oplus_e \det \R H(e)) \cong \det(\R E(X))\otimes \bigotimes_e
\det(\R H(e)).\eqno{(2)}$$
The second isomorphism follows since both expressions are determined by ordering and orienting all edges.

If $X$ is connected, $H_0(X;\mathbb R)\cong \R$ has a canonical (ordered!) basis.  Combining this
observation with isomorphisms (1) and (2) gives
	$$\det(H_1(X;\mathbb R))\otimes \det(\R V(X))\cong \det(\R E(X))\otimes \bigotimes_e \det(\R
H(e)).\eqno{(3)}$$

Now note that $\R E(X)$ and $\R V(X)$ have canonical unordered bases, which can be used to
identify
$\det(\R E(X))$ and
$\det(\R V(X))$ with their duals.  Since $V^*\otimes V$ is canonically isomorphic to $\R$, we can use this
fact to ``cancel" copies of $\det(\R E(X))$ or $\det (\R V(X))$, effectively moving them from one side of a
canonical isomorphism to the other.   In particular, from equation (3) we get the desired
canonical isomorphism
 $$\det(H_1(X;\mathbb R))\otimes \det(\R E(X))\cong \det(\R V(X))\otimes \bigotimes_e \det(\R
H(e)).\eqno{\qed}$$

Other equivalent notions of orientation, which we will use for particular types of graphs, are given in the next proposition 
and corollaries.  First we record an easy but useful lemma.

\begin{lemma}[Partition Lemma] \lbl{partition} Let $N$ be a finite set, and $P=\{P_1,\ldots,P_k\}$ a
partition of
$N$.  Then there is a canonical isomorphism 
$$\det(\R N)\cong \otimes \det(\R P_i)\otimes \det\big(\bigoplus_{|P_i| odd}\R\big),$$
which is independent of
the ordering of the $P_i$.
\end{lemma}

\begin{proof}  For $x_i$ in $N$, the map ``regroups" $x_1\wedge\ldots\wedge x_{|N|}$ so that all of the $x_i$
which are in $P_1$ come first, etc.  \end{proof}

\begin{proposition}\lbl{orient2} For $X$ a connected graph,  let $H(v)$ denote the set of
half-edges adjacent to a vertex $v$ of $X$. There is a canonical isomorphism 
$$
 \det(H_1(X;\mathbb R))\otimes \det(\R E(X)) \cong \bigotimes_{v\in V(X)}\det(\R
H(v))\otimes\det(\bigoplus_{|H(v)| even} \R) 
$$
\end{proposition}
\begin{proof}  
By   Lemma   ~\ref{partition},  grouping half-edges
according to the edges they form gives an isomorphism
$$\bigotimes_{e\in E(X)}\det(\R H(e))\cong\det\R H(X);$$
on the other hand,
 grouping according to the vertices to which they are adjacent gives an isomorphism 
$$\det\R H(X)\cong\bigotimes_v\det(\R H(v))\otimes \det(\bigoplus_{|H(v)| odd} \R v).$$

 Combining these isomorphisms, and substituting into
the canonical isomorphism of Proposition ~\ref{orient}, we get 
\begin{align*}
 \det(H_1(X;\mathbb R))&\otimes \det(\R E(X))
\\&\cong \det(\R V(X))\otimes\bigotimes_{v\in V(X)}\det(\R
H(v))\otimes\det(\bigoplus_{|H(v)| odd} \R v)\\ &\cong \bigotimes_{v\in V(X)}\det(\R
H(v))\otimes\det(\bigoplus_{|H(v)| even} \R v) 
\end{align*}
The last isomorphism follows from the fact that $$\det \R V(X)\cong \det(\bigoplus_{|H(v)| even} \R
v)\otimes \det(\bigoplus_{|H(v)| odd}\R v)$$ which, in turn, follows by the partition lemma combined
with the observation that a graph cannot have an odd number of vertices of odd valence and an odd
number of vertices of even valence.
\end{proof}
 
\begin{corollary} Let $T$ be a trivalent graph.  Then
an orientation on
$T$ is equivalent to a cyclic ordering of the edges incident to each vertex.
\end{corollary}

\begin{proof}
Since $T$ is trivalent, all vertices have odd valence, so the orientation is determined by
ordering the sets $H(v)$, up to cyclic (i.e.\ even) permutation, at each vertex $v$.
\end{proof}
This equivalence was mentioned in Kontsevich's papers \cite{K1},\cite{K2} and is also an
ingredient in the isomorphism between commutative graph \emph{cohomology} and the diagram
algebras arising in the study of finite type invariants of three-manifolds and knots.

\begin{corollary}\lbl{tree} Let $T$ be a connected binary tree.  An orientation on $T$ is equivalent to
an  ordering of the edges of $T$, or to a cyclic ordering of the
edges incident to  each interior vertex of $T$.
\end{corollary} 

\proof
Since $T$ is a connected binary tree, 
$H_1(T)$ is zero, and all of the vertices have odd valence (1 or 3), so that the isomorphism in the statement of  of Proposition
~\ref{orient2} reduces to:
$$\det(\R E(T))\cong \bigotimes_v\det(\R H(v))\eqno{(5)\quad \qed}$$

Note that a cyclic ordering of the edges incident to each interior vertex of a tree can be thought of as an
embedding of the tree into the plane.

\subsubsection{Chain groups}

We can now define the chain groups of the graph complex associated to a cyclic operad. 

\begin{definition} A vertex $v$ of a graph is 
\emph{$\O$-colored} if the half-edges incident to $v$ are identified with the
legs of an $\O$-spider. An \emph{$\O$-graph} is an oriented graph without univalent vertices   with an
$\O$-spider coloring each vertex.   
\end{definition} 
We represent an  $\O$-graph pictorially as in Figure~\ref{ograph}.
\begin{figure}[ht!]
\begin{center}
\includegraphics[height=3cm]{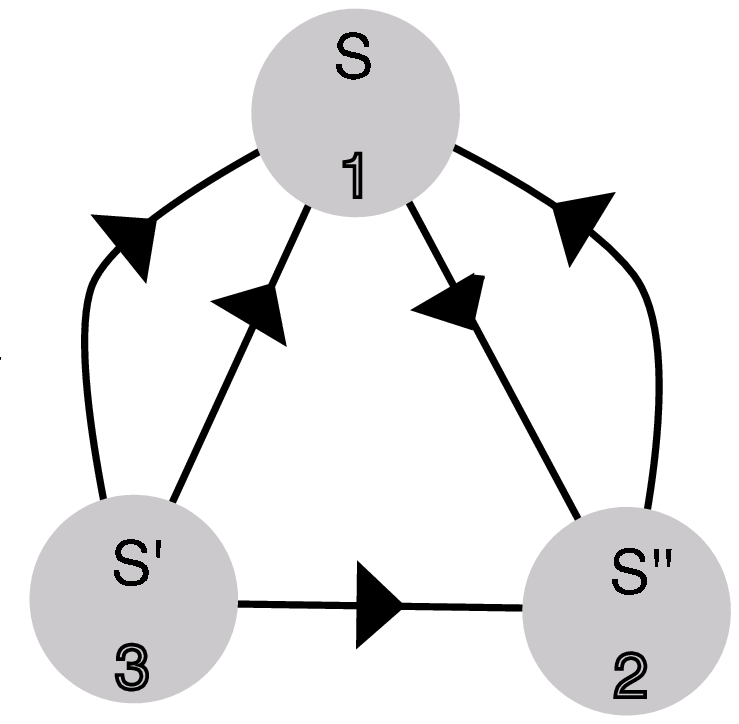}
\caption{$\O$-graph}\lbl{ograph}
\end{center}
\end{figure}

\begin{definition} The \emph{group of $k$-chains} $\G_k$ is a quotient of the vector space spanned
by
$\O$-graphs with $k$ vertices: 
$$\G_k= \mathbb R\{\text{$\O$-graphs with $k$ vertices}\}/relations$$
where the relations are of two kinds:
\begin{itemize}
\item[(1)] (Orientation)  $(\X,or)=-(\X,-or)$
\item[(2)] (Vertex linearity) If a vertex $v$ of $\X$ is colored by the element $S_v=aS+bT$, where
$a,b\in\mathbb R$ and $S,T\in \Sp[m]$, then $\X=a\X_S+b\X_T$, where $\X_S, \X_T$ are formed by coloring $v$ by $S$
and
$T$ respectively.
\end{itemize}
\end{definition}
Thus $\G_k$ is spanned by $\O$-graphs with $k$ vertices, each colored by a basis element of $\Sp$.  We set $$\G=\bigoplus_{k\geq 1}\G_k.$$

We  also define the \emph{reduced chain groups} $\overline{\G}_k$ to be $\G_k$ modulo the subspace of
graphs that have at least one vertex colored by $1_\O$.  
 In the commutative, associative and Lie cases,
 $\overline{\G}_k$ is spanned by $\O$-graphs without bivalent
vertices, since $\O[1]$ is spanned by $1_\O$ in these three cases.

\subsubsection{Hopf algebra structure}
  Both $\G$ and $\overline\G$ have a   Hopf algebra structure whose product is given by disjoint
union. More precisely, $X\cdot Y$ is defined to be the disjoint union $X\cup Y$ where the orientation
is given by shifting the vertex ordering of $Y$ to lie after that of $X$. The coalgebra structure is
defined so that connected graphs are primitive; the comultiplication is  then extended
multiplicatively to disjoint unions of graphs. The multiplicative unit is the empty graph, and the counit is dual to this unit. The
antipode reverses the orientation of a graph.  

The primitive parts (i.e.\ the subspaces of $\G$ and $\overline\G$ spanned by connected
graphs) will be denoted $P\G$ and $P\overline\G$ respectively.

\subsubsection{Boundary map}

Let $\X$ be an $\O$-graph, with underlying oriented graph $X$, and let $e$ be an edge of $X$.   We define a new
$\O$-graph
$\X_e$ as follows.  If $e$  is a loop, then $\X_e$ is zero.  If $e$ has distinct endpoints $v$ and $w$, then the underlying graph
$X_e$ is the graph obtained from
$X$ by collapsing
$e$. The orientation on
$X_e$ is  determined by the following rule: choose a representative of the orientation on $X$ so that $v$
is the first vertex, $w$ is the second vertex, and $e$ is oriented from $v$ to $w$.  The orientation
on
$X_e$ is then induced from that of $X$:  the uncollapsed edges are oriented as they were in $X$, the new
vertex resulting from collapsing $e$ is first in the vertex ordering, and the other vertices retain their
ordering. 
The $\O$-colorings at all the vertices besides $v$ and $w$ stay the same.
 Let $S_v$, $S_w$ be the $\O$-spiders coloring $v$ and $w$ respectively, with legs $\lambda$ of $S_v$ and $\mu$ of $S_w$ identified with the two
half-edges of $e$.
  Then the $\O$-spider $(S_v,\lambda)\circ(\mu,S_w)$ colors the vertex obtained by collapsing 
$e$ (Figure~\ref{collapse}).
\begin{figure}[ht!]\begin{center}
\includegraphics[height=3cm]{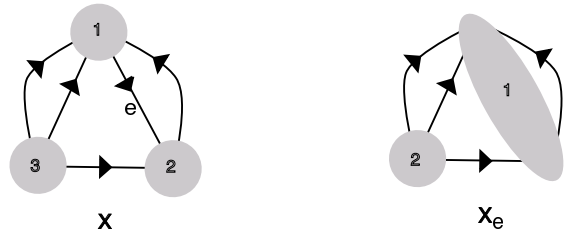}
\caption{Edge collapse}\lbl{collapse}
\end{center} 
\end{figure}

With this orientation convention, and using the associativity axiom of operads, the map
$$\bdry_E(\X) = \sum_{e\in X} \X_e$$ 
is a boundary operator. This makes
$\G$ into a chain complex, and we have

\begin{definition} The \emph{$\O$-graph homology} of the cyclic operad $\O$ is the homology of $\G=\oplus \G_k$ with respect to the boundary operator
$\bdry_E$.
\end{definition}
Note that  $\overline \G$, $P\G$ and $P\overline \G$ are all chain complexes with respect to
$\bdry_E$.   

We conclude this section with a nice observation, which we won't actually need, and whose proof is left
to the reader.

\begin{proposition}
The assignment $\O\mapsto \G$ respects morphisms, and hence is a functor from cyclic operads
to chain complexes.
\end{proposition}

\subsection{The Lie algebra of a cyclic operad and its homology} 

  In this section we associate  a sequence of  Lie algebras to any cyclic operad.  We show that each of
these Lie algebras contains a symplectic Lie algebra
 as a subalgebra, and that under certain finiteness assumptions the Lie algebra homology may be computed using the subcomplex of  symplectic invariants in
the exterior algebra.

\subsubsection{Symplectic Lie algebra of a cyclic operad} For each integer $n\geq 1$ we will define a functor \{Cyclic
Operads\}$\to$ \{Symplectic Lie Algebras\}, sending
$\O$ to $\A$.  We then take a limit as $n\to\infty$ to obtain an infinite-dimensional Lie algebra $\a\O_\infty$.
 
Fix a $2n$-dimensional real vector space $V_n$   basis
$B_n=\{p_1,\ldots,p_n,q_1,\ldots,q_n\}$ corresponding to the standard symplectic form $\omega$.
Given a cyclic operad $\O$, the idea is to form  a Lie algebra  by putting elements of $V_n$ on the legs of
$\O$-spiders, and defining the bracket of two such objects by summing over all possible matings, with
coefficients determined by the symplectic form.
\begin{figure}[ht!]
\begin{center}
\includegraphics[height=3cm]{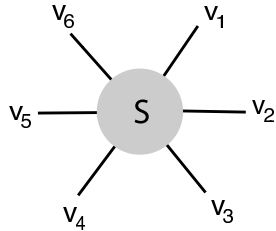}
\end{center}
\caption{Symplecto-spider}\lbl{symplecto}
\end{figure}

Formally, we ``put elements of $V_n$ on the legs  of an $\O$-spider" via a coinvariant
construction like the one used to define spiders, i.e.\ we 
 set
$$\A = \bigoplus_{m=2}^\infty (\Sp[m]\otimes V_n^{\otimes m})_{\Sigma_m},$$ 
where the symmetric group $\Sigma_m$ acts simultaneously on $\Sp[m]$ and on $V_n^{\otimes m}$.
We will refer to elements of $\A$ of the form $[S\otimes v_1\otimes\ldots\otimes v_m],$  where $S$ is an $\O$-spider, as
\emph{symplecto-spiders} (Figure~\ref{symplecto}).

The bracket of two  symplecto-spiders is defined as follows. Let $\ss_1=[S_1\otimes v_1\otimes\ldots\otimes v_m]$ and
$\ss_2=[S_2\otimes w_1\otimes\ldots\otimes w_l]$  be two symplecto-spiders.  Let $\lambda$ be a leg of $S_1$ and $\mu$
a leg of $S_2$, with associated elements $v_\lambda, w_\mu\in
 V_n$. Recalling that $\omega$ is the symplectic form on $V_n$, define
$(\ss_1,\lambda)!(\mu,
\ss_2)$ to be $\omega(v_\lambda,w_\mu)$ times the symplecto-spider obtained by  mating $S_1$ and $S_2$ using $\lambda$ and
$\mu$, erasing the elements
$v_\lambda$ and
$w_\mu$, and retaining  the elements of $V_n$ on the remaining legs (see
Figure~\ref{exclamation}).
\begin{figure}[ht!]
\begin{center}
\includegraphics[height=3cm]{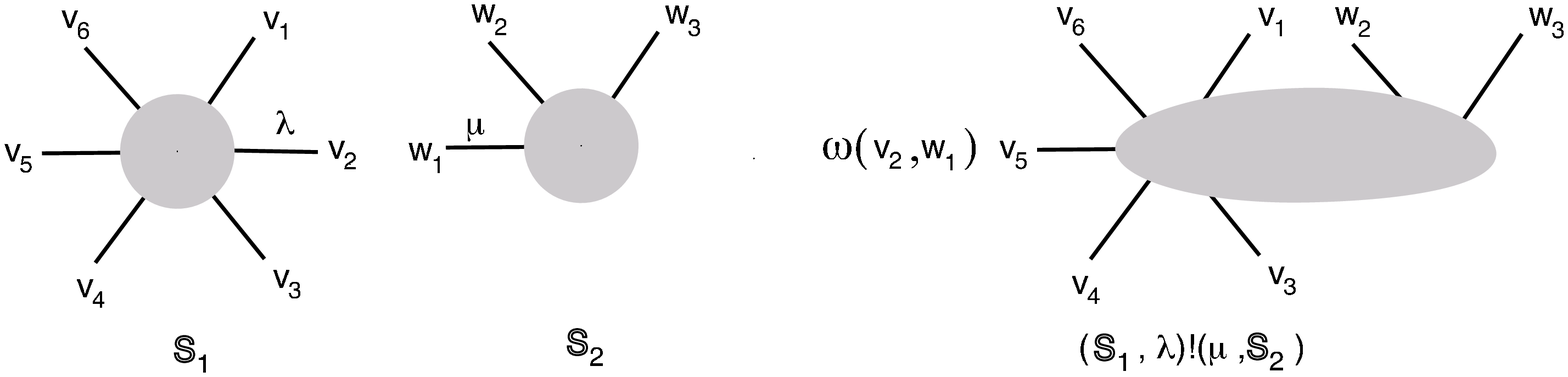}
\caption{$(\ss_1,\lambda)!(\mu, \ss_2)$}\lbl{exclamation}
\end{center}
\end{figure}

Now define the bracket   by setting
$$[\ss_1,\ss_2] = \sum_{\lambda\in \ss_1, \mu\in \ss_2} (\ss_1,\lambda)!(\mu, \ss_2); $$
and extending linearly to all of $\A$.

\begin{proposition}\lbl{antijacobi}
The bracket is antisymmetric and satisfies the Jacobi identity, for
any cyclic operad.
\end{proposition}

\begin{remark}\lbl{newremark}
A related construction appears in \cite{kapman}, Theorem 1.7.3, where, for any operad, a bracket on what might be called
the space of rooted spiders is defined. One sums over all ways of plugging the root (output)  of one spider into an input of another,
and then subtracts the results of doing this in the other order. In the cyclic case there is no specified root, so one would have to
sum over all choices of root; but then subtracting off the other order would be the same and give you a trivial operation.
The needed axiom in this case is that of an \emph{anticyclic} operad (\cite{GeKa, MSS}), which ensures that the order of
``plugging in" determines a sign. When $\O$ is a cyclic operad and $V$ is a symplectic vector space, the collection
$\O[m]\otimes V^{m+1}$ is an anticyclic operad. Thus the bracket defined here is a generalization to the anticyclic case
of the one defined by Kapranov and Manin.
\end{remark}

\begin{remark} 
This Lie algebra structure on $\A$ is quite natural.
Let  $T^{\geq 2}V_n$ denote the tensor algebra in degrees
$\geq 2$. $T^{\geq 2}V_n$ has a Lie bracket induced by the symplectic form,
 and we can give $\Sp$ an abelian Lie algebra structure.
Give the tensor product of associative algebras $\Sp\otimes T^{\geq 2}V_n$
the natural bracket which is a derivation in each variable, and which extends the brackets on each
tensor factor. 
Then the natural map $\Sp\otimes T^{\geq 2}V_n\to \A$ is a Lie algebra
homomorphism.
\end{remark}

 In the commutative case the Lie algebra $\A$ may be identified with the Lie algebra of  polynomials with no constant or
linear term in the variables $p_1\ldots p_n,q_1,\ldots,q_n$, under the standard Poisson bracket of
functions.

 In the Lie case the Lie algebra $\A$ is isomorphic to the Lie algebra $D_*(V_n)\otimes \R$  that has arisen in the study of the mapping class group
(see
\cite{Levine},
\cite{Morita}). Here $D_k(V_n)$ is
defined to be the kernel of the map $V_n\otimes L_{k}(V_n) \to L_{k+1}(V_n)$ sending $v\otimes x \mapsto [v,x]$, where
$L_*(V_n)$ refers to the free Lie algebra on $V_n$.

We record the functoriality of our construction without proof.
 
\begin{proposition}
For each $n\geq 1$, the assignment $\O\mapsto \A$  respects morphisms, and hence is a functor from cyclic operads
to Lie algebras.
\end{proposition}

Note that bracketing    a symplecto-spider with $2$ legs  and one with $m$ legs results in a sum of (at most two)
symplecto-spiders, each with
$m$ legs.  In particular, the subspace of
$\A$ spanned by symplecto-spiders with two legs forms a Lie subalgebra of $\A$.  If we consider only (two-legged)
symplecto-spiders with vertex colored by the identity $1_\O$, we obtain an even smaller subalgebra, denoted $\A^0$.
The next proposition identifies $\A^0$ with the symplectic Lie algebra $\sy(2n)$.

\begin{proposition}\lbl{symplectic}
Let $S_0$ denote the (unique) $\O$-spider colored by the identity element $1_\O$. 
The subspace $\A^0$ of $\A$ spanned by symplecto-spiders of the form $[S_0\otimes v \otimes w]$ is a Lie
subalgebra isomorphic to
$\sy (2n)$.
\end{proposition}
\begin{proof}
The map $\A^0\to S^2V$ sending $[S\otimes v\otimes w]$ to $vw$ is easily checked to be a Lie algebra
isomorphism, where the bracket on $S^2V$ is the Poisson bracket.  

Recall that $\sy(2n)$ is the  set of $2n\times 2n$ matrices $A$ satisfying $AJ+JA^T=0$, where 
$J= \begin{pmatrix} 0 & I  \\ -I &0\end{pmatrix}$.  The symmetric algebra $S^2V$ can be identified with the subspace of
$V\otimes V$ spanned by elements of the form $v\otimes w+ w\otimes v$.  Consider the composition of isomorphisms 
$$V\otimes V\to V^*\otimes V\to Hom(V,V)$$
where the first map is induced by the isomorphism $V\to V^*$ given by $v\mapsto \omega(v,\_)$. 
Tracing through these isomorphisms, we see that 
\begin{align*}
p_iq_j\leftrightarrow p_i\otimes q_j+q_j\otimes p_i &\mapsto \begin{pmatrix}  -E_{ji}&0\\ 0&E_{ij} \end{pmatrix}\\
p_ip_j\leftrightarrow  p_i\otimes p_j+p_j\otimes p_i &\mapsto \begin{pmatrix} 0&0\\E_{ij}+E_{ji}&0 \end{pmatrix}\\
  q_iq_j \leftrightarrow q_i\otimes q_j+q_j\otimes q_i &\mapsto \begin{pmatrix}  0&-E_{ij}-E_{ji}\\0&0 \end{pmatrix}
\end{align*}
where $E_{ij}$ is the $n\times n$ matrix with $(i,j)$-entry equal to 1 and zeroes elsewhere.  It is now straightforward to
check that this gives a Lie algebra isomorphism $S^2V\to \sy(2n)$.
\end{proof}

The subalgebra $\A^0\cong \sy(2n)$ acts on $\A$ via the bracket (i.e.\ the adjoint action). Using the remark following
Proposition~\ref{antijacobi}, we see that the
$\sy(2n)$ action on
$ (\Sp[m]\otimes V_n^{\otimes m})_{\Sigma_m}$ is given  by $\xi\cdot[S\otimes v_1\otimes\cdots\otimes v_m] =
\sum_{i=1}^m [S\otimes v_1\otimes\cdots\otimes (\xi\cdot v_i)\otimes\cdots\otimes v_m]$, where $V_n$ has the
standard $\sy (2n)$-module structure. 

The natural inclusion $V_n\to V_{n+1}$  induces an inclusion $\a\O_n\to \a\O_{n+1}$ of Lie algebras, which
is compatible with the inclusion $\sy(2n)\to \sy(2(n+1))$.
\begin{definition}
The infinite dimensional symplectic Lie algebra
$\a\O_\infty$ is the direct limit 
$$\a\O_\infty = \lim_{n\to\infty} \A.$$
\end{definition}

\subsubsection{Lie algebra homology}
The Lie algebra homology of $\A$ is computed
from the exterior algebra $\wedge  \A$ using the standard Lie boundary operator  $\bdry_n\colon
\wedge^k~\A \to
\wedge^{k-1} \A$ defined by
$$\bdry_n(\ss_1\wedge\ldots\wedge \ss_k)=\sum_{i<j} (-1)^{i+j+1}[\ss_i,\ss_j]\wedge
\ss_1\wedge\ldots\widehat{\ss_i}\wedge\ldots\wedge\widehat{\ss_j}\wedge\ldots\wedge \ss_k.$$
The map $\wedge\A\to \wedge\a\O_{n+1}$ induced by the natural inclusion is a chain map, 
so that 
$$H_k(\a\O_\infty;\R)=\lim_{n\to\infty} H_k(\A;\R).$$

\begin{proposition}
$H_k(\a\O_\infty;\R)$ has the structure of a Hopf algebra.
\end{proposition}
\begin{proof}
To define the product  $ H_*(\a\O_\infty)\otimes H_*(\a\O_\infty)\to H_*(\a\O_\infty)$, 
consider  maps $E\colon B_\infty\to B_\infty$ sending $p_i\mapsto p_{2i}$ (resp. $q_i\mapsto 
q_{2i}$), and $O \colon B_\infty\to B_\infty$ sending $p_i\mapsto p_{2i-1}$
(resp
$q_i\mapsto q_{2i-1}$).  These induce   maps $E$ and $O$ on $\a\O_\infty$.  The product   on
$H_*(\a\O_\infty)$ is induced by the map
$$\a\O_\infty\oplus\a\O_\infty\to\a\O_\infty$$
which sends $x\oplus y$ to $E(x)+O(y)$.

The coproduct $ H_*(\a\O_\infty)\to H_*(\a\O_\infty)\otimes H_*(\a\O_\infty) $ is induced by the
diagonal map
$\a\O_\infty\to\a\O_\infty\oplus\a\O_\infty$. 
More explicitly, the coproduct is induced by the  map $\wedge \a\O_\infty\to \wedge \a\O_\infty\otimes \wedge \a\O_\infty$
sending 
$$ \ss_1\wedge\ss_2\wedge\cdots\wedge \ss_k\mapsto\sum_{[k]=I\cup J}\epsilon(I,J) \ss_I\otimes
\ss_J,$$  where the sum is over unordered partitions of $[k]=\{1,\ldots, k\}$, where
$\ss_I=\ss_{i_1}\wedge\ldots\wedge
\ss_{i_{|I|}}$ if $I$ consists of $i_1<i_2<\ldots<\ldots i_{|I|}$, and where $\epsilon(I,J)$ is a sign determined
by the equation $\ss_1\wedge\ldots\wedge\ss_k = \epsilon(I,J)\ss_I\wedge\ss_J$. 

  The unit is $1\in \R\cong\wedge^0\a\O_\infty$, and the counit
is dual to this.  
\end{proof}

\subsubsection{The subcomplex of $\sy(2n)$-invariants} In the remainder of this paper, we assume that the
vector spaces $\O[m]$ are  finite-dimensional.  
In this section we show that in this case the homology of
$\a\O_n$ is computed by the subcomplex $(\wedge \A)^{\sy(2n)}$ of $\sy(2n)$ invariants (where an $\sy(2n)$
``invariant" is an element which is killed by every element of $\sy(2n)$.)  

In general, the exterior algebra $\wedge \A$ breaks up into a direct sum of 
pieces
$\Lambda_{k,m},$ spanned by wedges of $k$ symplecto-spiders with a total of $m$ legs:
\begin{align*}
\wedge&\A=\bigoplus_{k,m} \Lambda_{k,m}\\
&=\bigoplus_{k,m}(\bigoplus_{m_1+\ldots+m_k=mm+} (\Sp[m_1]\otimes V^{\otimes
m_1})_{S_{m_1}}\wedge\ldots\wedge (\Sp[m_k]\otimes V^{\otimes m_k})_{S_{m_k}}).
\end{align*}
If the vector spaces $\O[m]$ are finite-dimensional, then these pieces $\Lambda_{k,m}$ are all finite-dimensional, as are
the following subspaces:
\begin{definition}  The \emph{(k,m)-cycles} $Z_{k,m}$, the \emph{(k,m)-boundaries} $Z_{k,m}$ and the
\emph{(k,m)-homology} $H_{k,m}(\A;\R)$ are defined by
$$ 
Z_{k,m} =\Lambda_{k,m}\cap ker(\bdry_n),\quad
 B_{k,m} =\Lambda_{k,m}\cap im(\bdry_n),\quad
 H_{k,m}(\A;\R) =Z_{k,m}/B_{k,m}.
$$
\end{definition}

 With this definition, we have
$$  H_k(\A;\R)=\bigoplus_m H_{k,m}(\A;\R).$$ 

\begin{proposition}\lbl{invariants} 
The invariants $(\wedge \A)^{\sy(2n)}$ form a subcomplex of $\wedge\A$.
If $\O[m]$ is finite dimensional for all $m$, then the
inclusion $(\wedge \A)^{\sy(2n)}\to\wedge \A$ is an isomorphism on homology.
\end{proposition} 
\begin{proof}
The first statement follows since $\bdry_n$ is an $\sy(2n)$-module morphism.

The proof of the second statement depends on the following remark:  every
finite dimensional $\sy(2n)$-module, $E$, decomposes as a direct sum of $\sy(2n)$-modules in the
following way:
$$E= E^{\sy(2n)}\oplus
\sy(2n)\cdot E.$$
(Proof: $\sy(2n)$ is well-known to be reductive, which means that for every finite dimensional module $E$,
 every submodule $E^\prime \subset E$ has a complementary submodule $E^{\prime\prime}$ with
$E=E^\prime\oplus E^{\prime\prime}$. Thus $E=E^{\sy(2n)}\oplus E^{\prime\prime}$. Now $E^{\prime\prime}$
is a direct sum of simple modules on which $\sy(2n)$ acts nontrivially. Therefore, $\sy(2n)\cdot
(E^{\sy(2n)}\oplus E^{\prime\prime})=\sy(2n)\cdot E^{\prime\prime} = E^{\prime\prime}$, and the proof is
complete.)

By hypothesis  $\O[m]$ is finite-dimensional, so that $\is_{k,m}$ is a finite dimensional $\sy (2n)$-module. Since
$\sy(2n)$ is simple, this means that
$\is_{k,m}$ decomposes as a direct sum $\is_{k,m}=\is^{\sy(2n)}_{k,m}\oplus \sy(2n)\cdot \is_{k,m}$.
Since the boundary is an $\sy(2n)$ module morphism, the space of $(k,m)$-cycles $Z_{k,m}$ and the space of
$(k,m)$-boundaries $B_{k,m}$ are both $\sy(2n)$ modules. Therefore, they decompose as $
Z^{\sy(2n)}_{k,m}\oplus \sy(2n)\cdot Z_{k,m}$ and $ B^{\sy(2n)}_{k,m}\oplus \sy(2n)\cdot
B_{k,m}$ respectively. Thus the homology 
\begin{align*}
H_k(\A) &= \bigoplus_{m} H_{k,m}(\A)\\
        &= \bigoplus_{m} Z_{k,m}/B_{k,m}\\ 
        &= \bigoplus_m \frac{Z^{\sy(2n)}_{k,m}}{B^{\sy(2n)}_{k,m}}\oplus \frac{\sy(2n)\cdot Z_{k,m}}{\sy(2n)\cdot
                 B_{k,m}}
\end{align*}
Hence it suffices to show that $\sy(2n)\cdot Z_{k,m}=\sy(2n)\cdot
B_{k,m}$. For every $\xi\in \sy(2n)$ define $i^k_\xi: \wedge^{k-1}\A\to \wedge^k \A$ by $a\mapsto
\xi\wedge a$. Then one easily checks that for $a\in\wedge^k \A$,  $\xi\cdot a =
(\bdry_ni_\xi^{k+1}+i_\xi^k\bdry_n)(x)$. Thus, if $\xi\in\sy(2n)$ and $z\in Z_{k,m}$, then $\xi\cdot z =
\bdry_n i^{k+1}_\xi(z)$. Hence $\sy(2n)\cdot Z_{k,m} \subset B_{k,m}$, which implies
$[\sy(2n),\sy(2n)]\cdot Z_{k,m}\subset \sy(2n)\cdot B_{k,m}$. Since $\sy(2n)$ is simple, 
$[\sy(2n),\sy(2n)]=\sy(2n)$ and the proof is complete.
\end{proof}

\subsection{Relation between graph homology and Lie algebra homology}

Now we describe the construction at the heart of Kontsevich's proof,
namely the identification of $\sy(2n)$ invariants with oriented graphs.
There are three principle players that need to be introduced:
\begin{align*}
&\phi_n\colon\O\g\to \wedge \a\O_n\\
&\psi_n\colon\wedge \a\O_n \to \O\g\\
&M_n\colon\O\g\to\O\g
\end{align*}

\subsubsection{The map $\phi_n$}    We need to produce  wedges of symplecto-spiders from an $\O$-graph 
$\bf X$. Let $X$ be the underlying oriented graph, and fix an ordering of the vertices and directions on the edges representing the
orientation. A \emph{state}  of $\bf X$ is an assignment of an element of $B_n=\{p_1,\ldots,p_n,q_1\ldots,q_n\}$
to each half-edge of $X$ and a sign $\pm 1$ to each edge, subject to the following constraints:
\begin{itemize}
\item  if one half-edge of $e$ is labeled $p_i$, the other half-edge must
be labeled $q_i$, and vice versa;
\item
the sign on an edge is positive if the initial half-edge is labeled $p_i$, and negative if the initial half-edge is labeled $q_i$. 
\end{itemize}
Given a state $s$, The \emph{sign} of  $s$, denoted $\sigma(s)$, is the product of the signs on all edges, and we define an element  $\X\{s\}$ of
$\wedge\A$ as follows. Cut  each edge at its midpoint, thereby separating $\X$ into a disjoint union of symplecto-spiders  
$\ss_1,\ldots,\ss_k$ (where the subscript comes from the vertex ordering on $X$), and set $\X\{s\}= \ss_1\wedge\ldots \wedge\ss_k$  (See Figure
~\ref{state}). 

 Now define $\phi_n(\X)$ by summing over all possible states of $\X$:
$$\phi_n(\X)=\sum_s \sigma(s) \X\{s\}.$$
\begin{figure}[ht!]\begin{center}
\includegraphics[width=.9\hsize]{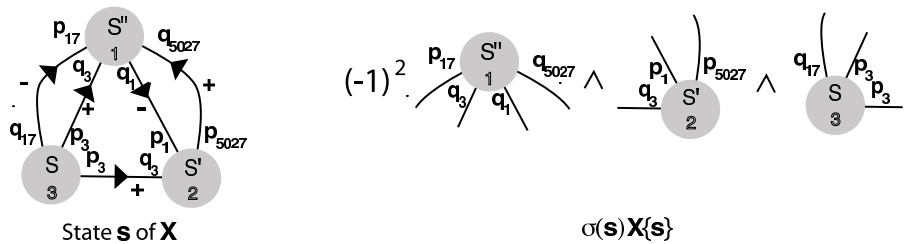} 
\caption{State of $\X$ and corresponding term of $\phi_n(\X)$}\lbl{state}
\end{center}
\end{figure}

\subsubsection{The map $\psi_n$}    We define $\psi_n$, too, as a state sum.  This time we need to
produce  $\O$-graphs from a wedge of symplecto-spiders
$\ss_1\wedge \ss_2\wedge \ldots\wedge \ss_k$. 
In this case a state  $\pi$ is a pairing of the legs of the spiders. We obtain a new $\O$-graph $(\ss_1\wedge\cdots\wedge
\ss_k)^\pi$ by gluing the spider legs together according to $\pi$, and orienting each edge arbitrarily.
Each edge of $(\ss_1\wedge\cdots\wedge \ss_k)^\pi$ carries an element $v_1\in V_n$ on its initial half-edge and $v_2\in V_n$ on its terminal half-edge.
We define the \emph{weight} of this edge to be $\omega(v_1,v_2)$, and denote the product of the weights of all edges by $w(\pi)$ (see Figure~\ref{pairing}).
\begin{figure}[ht!]\begin{center}
\includegraphics[width=.9\hsize]{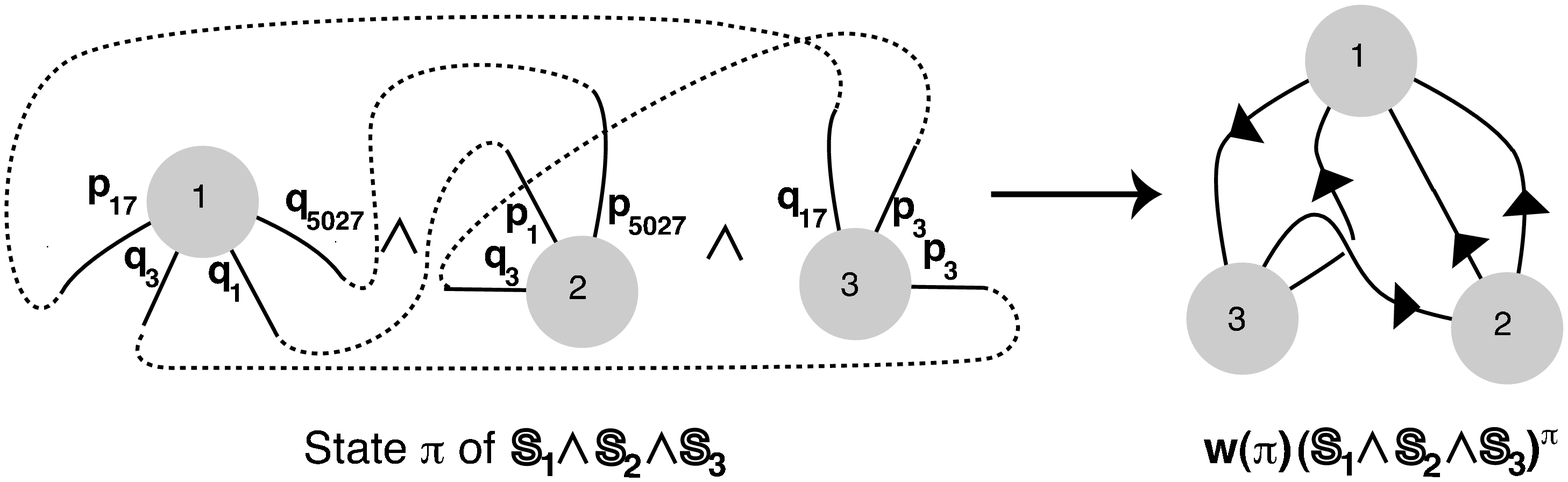}
\caption{A pairing $\pi$ and the resulting $\O$-graph}\lbl{pairing} 
\end{center}
\end{figure}
With this definition, the product $w(\pi) (\ss_1\wedge\ldots\wedge
\ss_k)^\pi$ is independent of the choice of edge-orientations, and we  define $$\psi_n(\ss_1\wedge \ss_2\wedge\ldots\wedge \ss_k)=\sum_\pi w(\pi)
(\ss_1\wedge\ldots\wedge
\ss_k)^\pi,$$
where the sum is over all possible pairings $\pi$.   Note that the vertex-linearity axiom in the
definition of $\G_k$ is required here to make $\psi_n$ linear.

\subsubsection{The map $M_n$}
As in the definition of $\psi_n$,
  a  pairing $\pi$ of the half-edges $H(X)$ of a graph $ X$ 
determines a new graph $X^\pi$, obtained by cutting all edges of $X$, then re-gluing the half-edges
according to
$\pi$.
 The \emph{ standard pairing} $\sigma$ pairs $x$ with $\bar x$ for all half-edges $x$ of
$X$; $X^\sigma$ is of course just $X$.  

If $X$ is oriented, there is an induced orientation on $X^\pi$, given as follows.
A pairing can be represented by a chord diagram on a set of vertices labeled by the half-edges of
$X$. The union of the chord diagrams for $\pi$ and for the standard pairing $\sigma$ forms a 
one-dimensional closed manifold $C(\pi)$, a union of circles.  Choose a representative for the orientation
of
$X$ so that the chords from the initial half-edge of each edge $e$ to the terminal half-edge of $e$ are
oriented coherently in each of these circles.  Now each edge of $X^\pi$ inherits a natural orientation
from each pair of half-edges determined by $\pi$.  The ordering of the vertices of $X^\pi$ is 
inherited from $X$.  If $X$ is the underlying graph of an $\O$-graph $\bf X$, we let $\X^\pi$ be the
induced $\O$-graph based on $X^\pi$ (see Figure~\ref{cpi}).   

Define $c(\pi)$ to be the number of components of $C(\pi)$.   and define the  map $M_n\colon \G\to \G$ by
$$M_n(\X)=\sum_{  \pi} (2n)^{c(\pi)}\X^\pi,$$ 
where the sum is over all possible pairings $\pi$ of $H(X)$.
\begin{figure}[ht!]\begin{center}
\includegraphics[width=.9\hsize]{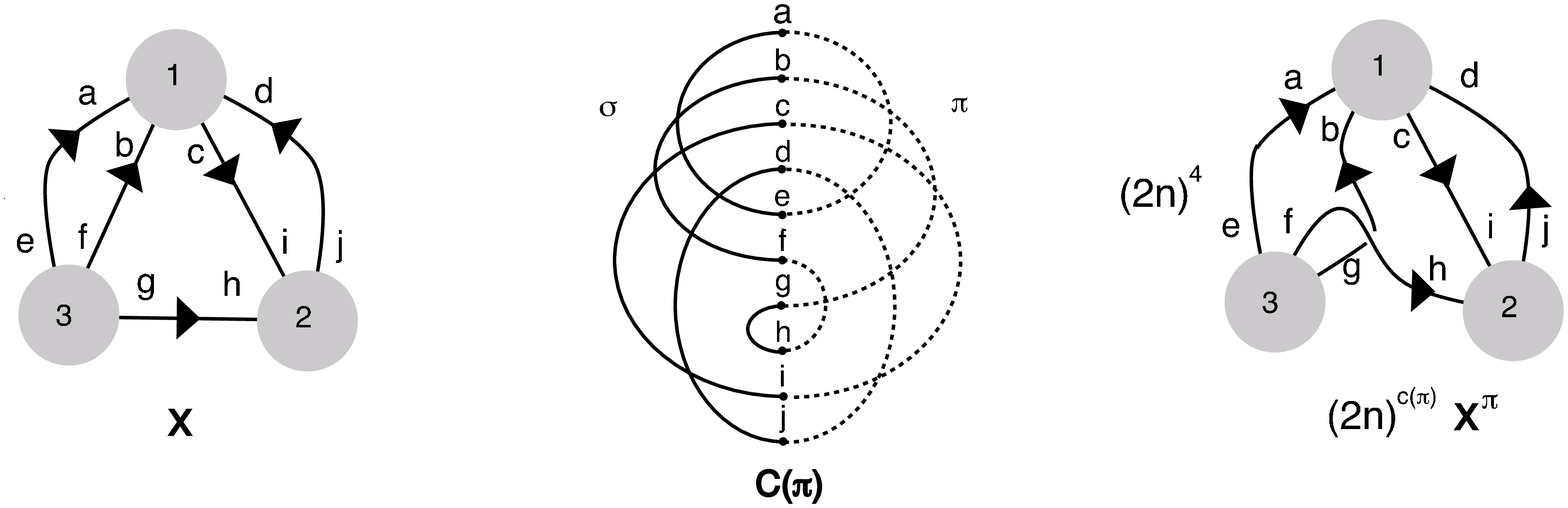} 
\caption{A term in $M_n(\X)$}\lbl{cpi}
\end{center}
\end{figure}
 
The map $M_n$ decomposes as a direct sum as follows.  Write $$\O\g=\bigoplus_{k,m} \G_{k,m}$$ where $\G_{k,m}$ is spanned by $\O$-graphs with  $k$
vertices and
$m$  half-edges (i.e $m/2$ edges). The  subspaces $\G_{k,m}$ are invariant under
$M_n$, and we denote by $M^{k,m}_n$ the restriction of $M_n$ to $\G_{k,m}$. 

\begin{proposition}\lbl{Mkl} If $\O[m]$ is finite-dimensional, then for large enough $n$, the restriction
$$M^{k,m}_n\colon \O\g_{k,m}\to\O\g_{k,m}$$
of $M_n$ is an isomorphism.  
\end{proposition}
\begin{proof}
If $\O[m]$ is finite-dimensional, then so is $\O\g_{k,m}$, and we can think of the restriction  $M^{k,m}_n$ of $M_n$ to
$\O\g_{k,m}$, as a matrix.  The matrix entries are polynomials in
$n$. The maximum of $c(\pi)$ occurs when $\pi$ is the standard pairing $\sigma$,
and is equal to $m$. Thus the diagonal entries are of the form $(2n)^m+${\it (lower degree terms)}, whereas the off diagonal
terms are all of lower degree. Therefore, for large enough $n$, the matrix is invertible, i.e.\ the map
$M^{k,m}_n$ is an isomorphism. 
\end{proof}

The three maps are related by the following Proposition.
\begin{proposition}\lbl{Mn}
$\psi_n\circ\phi_n = M_n$
\end{proposition}

\begin{proof}
Applying $\phi_n$ to an $\O$-graph $\bf X$ means that we are assigning elements of $B_n$ to the endpoints of the chord diagram for the standard pairing
$\sigma$ of $H(X)$, as on the left of Figure~\ref{proofpic}. 
\begin{figure}[ht!]\begin{center}
\includegraphics[height=6cm]{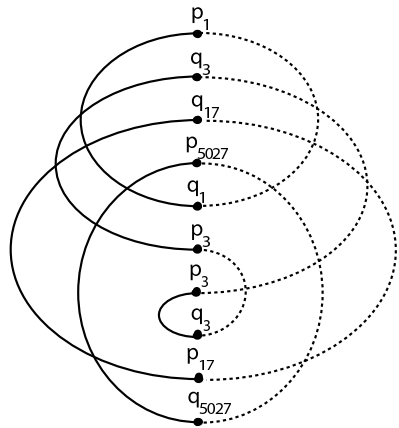}
\caption{Proof of Proposition~\ref{Mn}}\lbl{proofpic}
\end{center}
\end{figure}
Each chord must connect a pair $\{p_i, q_i\}$, for some $i$.
To then apply $\psi_n$, we consider all possible pairings $\pi$ of $H(X)$, and reglue to get $\O$-graphs $\X^\pi$. The weight $w(\pi)$ will
only be non-zero if every chord for $\pi$ also connects a pair $\{p_i,q_i\}$ for some $i$, as on the right of Figure~\ref{proofpic}. 
Thus the label of a vertex in $C(\pi)$ determines the labels of all the vertices of each connected component of $C(\pi)$: they must alternate
$p_i,q_i,p_i,q_i,...$ as you travel around a circuit. There are $2n=|B_n|$ choices of label for each component, so there are
$(2n)^{c(\pi)}$ possible pairings $\pi$ with non-zero weight.  Keeping track of the orientations, we see that each of these terms has
weight
$1$, so that 
the composition $\psi_n\circ\phi_n$ is exactly  equal to $M_n$. \end{proof}

  Recall that $\wedge\A=\bigoplus \is_{k,m}$, where $\Lambda_{k,m}$ is spanned by wedges of $k$
symplecto-spiders with a total of $m$ legs. Let
$\psi_{k,m}$ denote the restriction of
$\psi$ to
$\is_{k,m}$, and notice that
$\psi_{k,m}(\is_{k,m})\subset \O\g_{k,m}$. Similarly denote the restriction of $\phi_n$ to $\is_{k,m}$
by $\phi_{k,m}$ and notice that $\phi_{k,m}(\is_{k,m})\subset \O\g_{k,m}$. 
\begin{corollary}\lbl{psionto}
For $n$ sufficiently large with respect to fixed $k$ and $m$, 
\begin{itemize}
\item[\rm i)] the map $\psi_{k,m}\colon \is_{k,m}\to \O\g_{k,m}$ is onto, and
\item[\rm ii)] the map $\phi_{k,m}\colon \O\g_{k,m}\to \is_{k,m}$ is injective.
\end{itemize}
\end{corollary}

\subsubsection{Graphs and invariants}  The following proposition shows that the map $\phi_n$ gives instructions for constructing
an $\sy(2n)$-invariant from an $\O$-graph, and that all $\sy(2n)$-invariants are constructed in this way.

\begin{proposition}\lbl{imphi}
$ im(\phi_n)=(\wedge \A)^{\sy(2n)} $
\end{proposition}
\begin{proof}
In order to determine the $\sy(2n)$-invariants in $\wedge \A$, we first lift to the tensor algebra $T(\A)$.  
The quotient map $p\colon T(\A)\to \wedge \A$ sending $\ss_1\otimes\ldots\otimes \ss_k\mapsto \ss_1\wedge\ldots\wedge \ss_k$ is an $\sy(2n)$-module
map, so restricts to $p\colon (T(\A))^{\sy(2n)}\to (\wedge \A)^{\sy(2n)}.$  To see that $p$ is surjective, note that composition
$$\wedge \A   \stackrel{i}{ \to}   T \A \stackrel{p}{ \to}  \wedge \A, $$ 
where the map $i$ is defined by
$i(\ss_1\wedge\ldots\wedge \ss_k)= \frac{1}{k!}\sum_{\sigma\in \Sigma_m}
sgn(\sigma) \ss_{\sigma(1)}\otimes\ldots\ldots \otimes\ss_{\sigma(k)}$,  is the
identity. Since $i$ is also an $\sy(2n)$-module homomorphism,  
we get restrictions
$$ (\wedge \A)^{\sy(2n)}\stackrel{i}{ \to}(T(\A))^{\sy (2n)} \stackrel{p}{ \to}(\wedge \A)^{\sy(2n)}$$
whose composition is the identity. In particular the restriction of $p$ is onto.

We next lift even further.  Recall that $\A = \oplus_{m=2}^\infty(\Sp[m]\otimes V^{\otimes
m})_{\Sigma_m}$, where $V=V_n$.  Define $\widehat{\a\O}_n$ to be  $\oplus_{m=2}^\infty(\Sp[m]\otimes V^{\otimes
m})$.
The map $q\colon T(\widehat{\a\O}_n)\to T(\A)$ induced by the quotient maps $(\Sp[m]\otimes V^{\otimes m})
\to (\Sp[m]\otimes V^{\otimes m})_{\Sigma_m}$ is an $\sy (2n)$ module map, so restricts to $q\colon
T(\widehat{\a\O}_n)^{\sy(2n)}\to T(\A)^{\sy(2n)}$.  To see that this is surjective, consider the  composition
$$T(\A) \stackrel{j}{ \to} T(\widehat{\a\O}_n)\stackrel{q}{ \to} T(\A),$$
where $j$ is induced by the maps
$$(\Sp[m]\otimes V^{\otimes m})_{\Sigma_m}\to \Sp[m]\otimes V^{\otimes m}$$
sending $x$ to
$\frac{1}{m!}\sum_{\sigma\in \Sigma_m}\sigma\cdot x$. These maps are  well defined (they
send $x$ and $\sigma\cdot x$ to the same element) and the composition $q\circ j$ is the identity.
In addition, $j$ is  an $\sy(2n)$-module homomorphism, so
the restriction
$$(T(\A))^{\sy(2n)} \stackrel{j}{ \to} (T(\widehat{\A}))^{\sy(2n)}\stackrel{q}{ \to} (T(\A))^{\sy(2n)}$$  
is the identity; in particular, the restriction of $q$ is onto. 

We now compute the invariants in $T(\widehat{\a\O}_n)$. We have
$$T(\widehat{\a\O}_n)=\bigoplus_{k\geq 2,m\geq 1}\widehat{T}_{k,m},$$
where
$$\widehat{T}_{k,m} = \bigoplus_{m_1+\cdots+m_k=m} (\Sp[m_1]\otimes V^{\otimes
m_1})\otimes \cdots\otimes (\Sp[m_k]\otimes V^{\otimes
m_k}) $$
 The
action of $\sy (2n)$ on each $\Sp[k]\otimes V^{\otimes k}$ affects only the $V^{\otimes k}$ factor, so that

$(\widehat{T}_{k,m})^{\sy(2n)}$
$$= \bigoplus_{m_1+\cdots+m_k=m} \left((
V^{\otimes m_1}\otimes \cdots\otimes  V^{\otimes
m_k})^{\sy(2n)}\otimes (\Sp[m_1]\otimes\cdots\otimes\Sp[m_k]) \right) $$
Computation of the invariants $(V^{\otimes N})^{\sy(2n)}$ is a classical result of Weyl \cite{Weyl}.  A convenient 
way of parameterizing these is by
oriented  chord diagrams. 
Consider a set of $N$ vertices representing the tensor factors of $V^{\otimes N}$.
An oriented chord diagram is a choice of directed edges  that pair up the vertices; such a diagram
gives rise to an invariant as follows. 
 As in the definition of $\phi_n$, define a \emph{state} to be an assignment of pairs $\{p_i,q_i\}$ to the vertices joined by each chord.
The invariant is then the sum, over all possible signed states, of the associated elements of $V^{\otimes N}$.  
For example, when $N=4$, consider the
chord diagram $1\to 3, 4\to 2$. One state for this diagram assigns $p_1$ to the first vertex, $q_1$ to the third vertex, $q_6$ to the fourth vertex and
$p_6$ to the second vertex.  The corresponding term of $V^{\otimes 4}$  is 
$-p_1\otimes p_6\otimes q_1\otimes q_6$. The chord $4\to 2$  contributes a minus sign,   since the $p_6$ occurs at the head instead
of the tail of the chord.

Having identified the invariants in $T(\widehat{\a\O}_n)$, we follow the surjective maps 
$$
(T(\widehat{\A}))^{\sy(2n)}\to (T(\A))^{\sy(2n)}\to
(\wedge\A)^{\sy(2n)}$$
to identify the invariants in $\wedge \A$. 
Let $c\otimes(S_1\otimes\ldots\otimes S_k)\in(
V_n^{\otimes m_1}\otimes \cdots\otimes  V_n^{\otimes
m_k})^{\sy(2n)}\otimes (\Sp[m_1]\otimes\cdots\otimes\Sp[m_k])$
where $c$ is a chord diagram invariant.
Mapping to $$(\Sp[m_1]\otimes V_n^{\otimes
m_1})_{\Sigma_{m_1}}\otimes \cdots\otimes (\Sp[m_k]\otimes V_n^{\otimes
m_k})_{\Sigma_{m_k}}$$ has the effect of putting the vectors $p_i, q_i$  in the terms of $c$
on the legs of the spiders $S_1,\ldots, S_k$. 
The chord diagram $c$ gives instructions for gluing the spider legs to form an $\O$-graph $\X$, and each state
of $c$ induces a state of $\X$.  Summing over all states of $c$, we obtain the sum over all states of $\X$.
Mapping then to $\wedge\A$ explains the antisymmetry
of the vertex ordering of $\X$. 

The invariant just described is precisely $\phi_n(\X)$.
 Since these invariants span the space of all invariants, the proposition is
proved.
\end{proof}

\begin{corollary}\lbl{invariantiso}
For large enough $n$ with respect to fixed $k,m$, the restriction $$\psi_{k,m}\colon(\Lambda_{k,m})^{\sy(2n)}\to
\G_{k,m}$$ is an isomorphism.
\end{corollary}
\begin{proof}  This is immediate from Propositions \ref{Mkl}, \ref{Mn} and \ref{imphi}.
\end{proof}

\subsubsection{Isomorphism with  graph homology}

The map $\phi_n$ is \emph{not} a chain map, so does not induce an isomorphism on homology; this
was an oversight in Kontsevich's original papers \cite{K1},\cite{K2}.  We show instead that $\psi_n$ is a chain
map, and induces an isomorphism on homology after stabilization with respect to $n$.

\begin{proposition}\lbl{comm}
$\psi_n\colon\Lambda\A\to \O\g$ is a chain map.
\end{proposition}
\begin{proof}

We need to show the following diagram commutes. 
\begin{equation*}
\begin{CD} 
\wedge\A@>\psi_n>>\O\g  \\ 
@V{\bdry_n}VV @VV{\bdry_E}V \\ 
\wedge\A@>\psi_n>>\O\g
\end{CD}
\end{equation*}
Let $\ss_1\wedge\ldots\wedge \ss_k\in\wedge^k\A$. Then  
{\small\begin{align*}
&\psi_n\bdry_n(\ss_1\wedge\ldots\wedge \ss_k)\\
&\quad=\psi_n\left(\sum_{i<j}(-1)^{i+j+1}[\ss_i,\ss_j]\wedge
\ss_1\wedge\ldots\wedge\hat{\ss}_i\wedge\ldots
\wedge\hat{\ss}_j\wedge\ldots\wedge \ss_k \right)\\  
&\quad=\psi_n\left(\sum_{i<j}(-1)^{i+j+1}\sum_{\lambda\in \ss_i,\mu\in \ss_j}(\ss_i,\lambda)!(\mu,\ss_j)\wedge
\ss_1\wedge\ldots\wedge\hat{\ss}_i\wedge\ldots
\wedge\hat{\ss}_j\wedge\ldots\wedge \ss_k\right)\\
&\quad=\sum_{i<j}\sum_{\lambda\in \ss_i,\mu\in \ss_j}\\ 
&\hspace{.5in}\sum_{\pi}(-1)^{i+j+1} \,w(\pi)
\left((\ss_i,\lambda)!(\mu,\ss_j)\wedge
\ss_1\wedge\ldots\wedge\hat{\ss}_i\wedge\ldots
\wedge\hat{\ss}_j\wedge\ldots\wedge \ss_k\right)^\pi\tag{*}
\end{align*}}%
where $\pi$ runs over all pairings of legs of the $\ss_i$ other than $\lambda$ and $\mu$.

On the other hand
\begin{align*}
\bdry_E\psi_n(\ss_1\wedge\ldots\wedge \ss_k) &= \bdry_E\sum_\tau w(\tau)(\ss_1\wedge\ldots\wedge\ss_k)^\tau\\ 
&=\sum_\tau w(\tau)\sum_{e\in (\ss_1\wedge\ldots\wedge \ss_k)^\tau}(\ss_1\wedge\ldots\wedge\ss_k)^\tau_e\tag{1}\\ 
\end{align*}
 where $\tau$ runs over all pairings of the legs of the $\ss_i$.  Two legs $\lambda$ and $\mu$ form an edge $e$ of $(S_1\wedge\ldots\wedge
S_k)^\tau$ if and only if they are paired by
$\tau$, i.e.\
$\{\lambda,\mu\}\in\tau$. If $\lambda$ and $\mu$ are in the same $S_i$ for some $i$, then the edge $e$ is a loop, and
$(S_1\wedge\ldots\wedge S_k)^\tau_e=0$.  Therefore we can rewrite (1) as
\begin{align*}
&=\sum_\tau w(\tau)\sum_{\{\lambda,\mu\}\in  \tau}(\ss_1\wedge\ldots\wedge \ss_k)^\tau_{\lambda\cup \mu} \\
&=\sum_\tau \sum_{i<j}\sum_{\lambda\in \ss_i,\mu\in\ss_j, \{\lambda,\mu\}\in\tau}w(\tau)(\ss_1\wedge\ldots\wedge \ss_k)^\tau_{\lambda\cup
\mu}. \tag{2}
\end{align*}
Now set $\tau'=\tau-\{\lambda,\mu\}$, i.e.\ $\tau'$ is the pairing on the legs other than $\lambda$ and $\mu$ induced by $\tau$. Then
$w(\tau)=w(\tau')\omega(v_\lambda,v_\mu)$, and 
{\small\begin{align*}
w(\tau)(\ss_1\wedge&\ldots\wedge \ss_k)^\tau_{\lambda\cup
\mu}\\
&=(-1)^{i+j+1}w(\tau')\omega(v_\lambda,v_\mu)(\ss_i\wedge \ss_j\wedge \ss_1\wedge\ldots\wedge \hat
\ss_i\wedge\ldots\wedge\hat\ss_j\wedge\ldots\wedge
\ss_k)^\tau_{\lambda\cup\mu}\\
&=(-1)^{i+j+1}w(\tau') ((\ss_i,\lambda)!(\mu, \ss_j)\wedge \ss_1\wedge\ldots\wedge \hat \ss_i\wedge\ldots\wedge\hat\ss_j\wedge\ldots\wedge
\ss_k)^{\tau'}.
\end{align*}}%
Substituting into the sum (2) gives
{\small$$
=\sum_\tau \sum_{i<j}\sum_{\lambda\in \ss_i,\mu\in\ss_j,
\{\lambda,\mu\}\in\tau}(-1)^{i+j+1}\,w(\tau')((\ss_i,\lambda)!(\mu, \ss_j)\wedge \ss_1\wedge\ldots\wedge \hat
\ss_i\wedge\ldots\wedge\hat\ss_j\wedge\ldots\wedge
\ss_k)^{\tau'}
$$}%
Changing the order of summation now recovers the formula (*) for $\psi_n\bdry_n(\ss_1\wedge\ldots\wedge \ss_k)$.
\end{proof}
%

We now want to stabilize by letting $n\to\infty.$  We need

\begin{proposition}\lbl{stable} The system of maps $\psi_n$ gives rise to a map $\psi_\infty\colon\wedge\a\O_\infty\to \G$, i.e.\  
 the following diagram commutes:
\begin{equation*}
\begin{CD}
\wedge\A@>{\psi_n}>>\O\g\\
@VVV @VV=V\\
\wedge \a\O_{n+1}@>{\psi_{n+1}}>>\O\g
\end{CD}
\end{equation*}
\end{proposition}

\begin{proof}
Write $V_{n+1}=V_n\oplus V_1$, where $V_1$ is spanned by $p_{n+1},q_{n+1}$.
Commutativity of the diagram follows since the gluing instructions given by $\psi_{n+1}$ are the same as those of
$\psi_n$ if the symplecto-spiders happen to have legs labeled only by elements of $V_n$.
\end{proof}

We are now ready to prove the main theorem, identifying the homology of the infinite-dimensional Lie algebra with $\O$-graph homology.

\begin{theorem}\lbl{hopfiso}
The  map $\psi_\infty:\wedge\a\O_{\infty}\to \O\g$ induces a Hopf algebra isomorphism on homology
$H_*(\a\O_\infty ;\mathbb R)\cong H_*(\O\g)$.
\end{theorem}

\begin{proof} We have
$$H_*(\wedge\a\O_\infty;\mathbb R)=\lim_{n\to\infty}H_*(\wedge\a\O_n;\mathbb R)=\bigoplus_{k,m}\lim_{n\to\infty} H_{k,m}(\wedge\A;\mathbb R).$$ 

By Proposition~\ref{invariants}, $H_{k,m}(\wedge\A;\mathbb R)\cong H_{k,m}((\wedge\A)^{\sy (2n)}
;\mathbb R).$ By Corollary~\ref{invariantiso}, there is an $N$ such that
$(\is_{k,m})^{\sy(2n)},(\is_{k+1,m+2})^{\sy(2n)},$ and $ 
(\is_{k-1,m-2})^{\sy(2n)}$ are isomorphic via $\psi_n$ to $$\O\g_{k,m},\O\g_{k+1,m-2},\O\g_{k-1,m+2}$$ respectively
for $n>N$. The fact that these isomorphisms
respect the boundary map (Proposition~\ref{comm}) 
implies that  for $n>N$ $(\psi_n)_*\colon H_{k,m}(\wedge\A;\mathbb R)\to H_{k,m}(\O\g)$ is an isomorphism.  Since $H_{k,m}(\G)$ is independent of $n$, we
have 
$$\bigoplus_{k,m}\lim_{n\to\infty} H_{k,m}(\wedge\A;\mathbb R)
\cong\bigoplus_{k,m} H_{k,m}(\O\g)=H_*(\G), 
$$
showing that $\psi_\infty$ induces an isomorphism on homology.
%

To see that $\psi_\infty$ is a Hopf algebra isomorphism,   it suffices to show that the 
map on the chain level is a Hopf algebra morphism.  We first show that it is
compatible with the product structures, i.e.\ the following diagram commutes.
\begin{equation*}
\begin{CD}
\wedge(\a\O_\infty)\otimes\wedge(\a\O_\infty)@>>>\wedge\a\O_\infty\\
@VV{\psi_\infty\otimes\psi_\infty}V @VV{\psi_\infty}V\\
\G\otimes\G@>>>\G
\end{CD}
\end{equation*}
Let $x\otimes y\in\wedge(\a\O_\infty)\otimes\wedge(\a\O_\infty)$, where $x$ and $y$ are each wedges of
symplecto-spiders with legs labeled by $B_\infty=\cup B_n$  (i.e.\ the legs are labeled by the basis vectors
$\{p_i,q_i\}_{i=1}^{\infty}$).
 The top map sends $x\otimes y$ to $E(x)\wedge O(y)$. The
labels of
$E(x)$ are disjoint from those of $O(y)$, hence gluing up the legs with $\psi_\infty$ will act separately on $E(x)$ and
$O(y)$ yielding the disjoint union $\psi_\infty(E(x))\cup\psi_\infty(O(y))$. This says exactly that the diagram commutes.

Now we show it is a coalgebra morphism on the chain level:
\begin{equation*}
\begin{CD}
\wedge\a\O_\infty@>>> \wedge\a\O_\infty\otimes \wedge\a\O_\infty\\
@VV{\psi_\infty}V @VV{\psi_\infty\otimes\psi_\infty}V\\
\G @>>> \G\otimes\G
\end{CD}
\end{equation*}
Consider a wedge of $B_\infty$-labeled symplecto-spiders $\ss_1\wedge\ldots\wedge \ss_k$. 
Both routes around the rectangle involve partitioning the spiders and gluing up the legs. 
Let $\pi$ be a gluing, i.e.\ a pairing of the legs of the symplecto-spiders. Let $P=(I,J)$ be a partition of
the spiders into two sets. In order for $P$ to give a summand of $\Delta(\ss_1\wedge\ldots\wedge \ss_k)^\pi$,
$\pi$ should only pair legs of symplecto-spiders in the same part of the partition. In this case we say that
$\pi$ and $P$ are \emph{consistent}. 
On the other hand, a typical term of $\Delta(\ss_1\wedge\ldots\wedge \ss_k)$ is $\epsilon(I,J) \ss_I\otimes
\ss_J$. In order for $\pi$ to give a term of $(\psi\otimes\psi)(\epsilon(I,J) \ss_I\otimes \ss_J)$,
$\pi$ and $P$ should be consistent.
Therefore, going both ways around the rectangle, the summands are parameterized by consistent pairs $(P,\pi)$.

It remains only to check the coefficients.  Given a consistent pair $(P,\pi)$, we can subdivide $\pi$ into two pairings $\pi_I$ and $\pi_J$, which pair
the legs in the
$I$ and
$J$ symplecto-spiders respectively.
 Now the $(P,\pi)$ term in  $\psi_\infty\Delta(\ss_1\wedge\ldots\wedge \ss_k)$ is 
$\epsilon(I,J)(\ss_I\otimes \ss_J)^\pi=\epsilon(I,J)\ss_I^{\pi_I}\otimes \ss_J^{\pi_J}$.
The $(P,\pi)$ term in $\Delta\psi_\infty(\ss_1\wedge\ldots\wedge \ss_k)=\epsilon(I,J)\Delta\psi_\infty (\ss_I\wedge
\ss_J)$ is the $P$ term in $\epsilon(I,J)(\ss_I\wedge \ss_J)^\pi=\epsilon(I,J)\ss_I^{\pi_I}\wedge \ss_J^{\pi_J}$,
which is just $\epsilon(I,J)\ss_I^{\pi_I}\otimes \ss_J^{\pi_J}$.

We illustrate the commutativity of the coalgebra diagram  in Figure~\ref{coalgebra}.
\begin{figure}[ht!]\begin{center}
\includegraphics[width=\hsize]{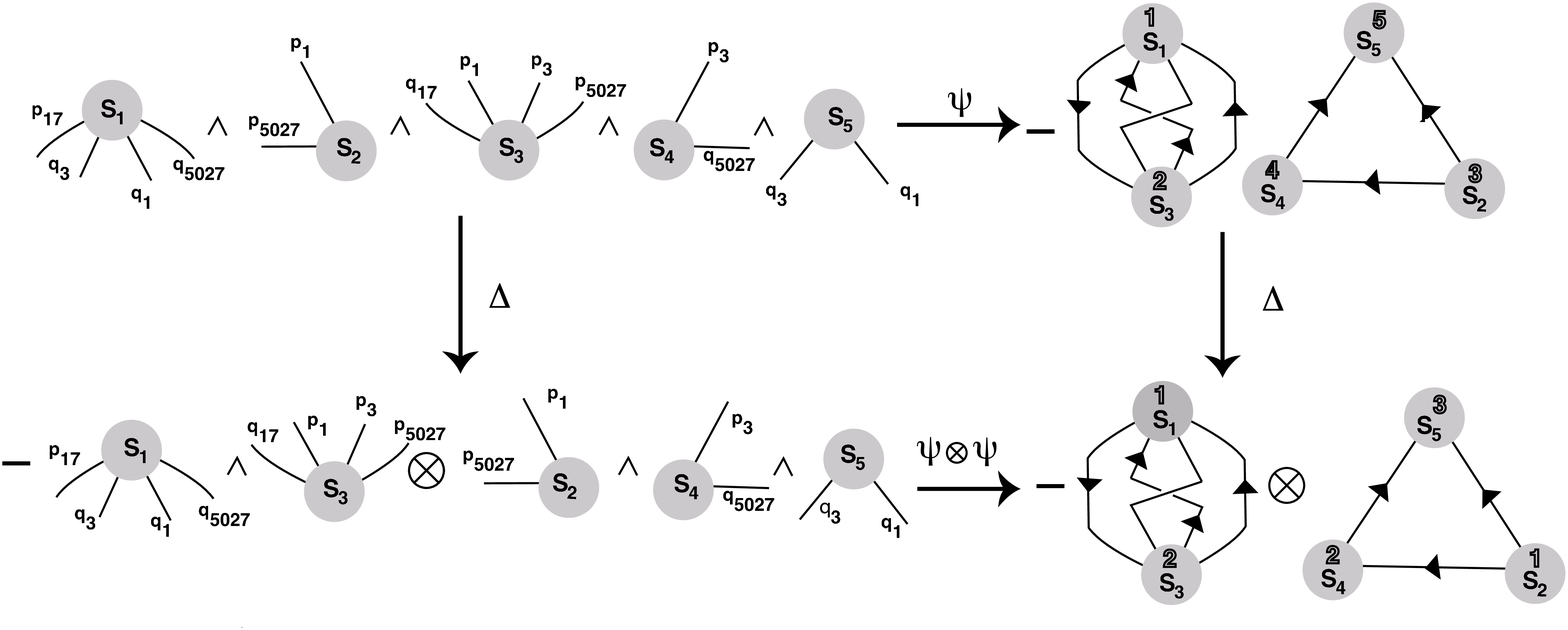}
\caption{Commutativity of the coalgebra diagram}\lbl{coalgebra}
\end{center}
\end{figure} 
In the upper right is a possible gluing, $\pi$, of the legs, with numbering of the vertices depicted by the
outlined numbers. We have switched the numbers $2$ and $3$, resulting in a minus sign. In the lower left is a
possible partition, $P$, of the symplecto-spiders, which is consistent with the gluing. In the lower right, we
have the result of applying $\pi$ to the lower left, which is the same as applying the partition $P$ to the
upper right.
\end{proof}

Recall that the primitive elements of $\G$ are the connected $\O$-graphs, spanning the subcomplex $P\G$.  If we denote the primitive elements of
$H_*(\a\O_\infty;\mathbb R)$ by $Prim\left(H_*(\a\O_\infty;\mathbb R)\right)$, we have

\begin{corollary}
$Prim\left(H_*(\a\O_\infty;\mathbb R)\right)\cong H_*(P\G)$.
\end{corollary}

\subsubsection{The primitive irreducible subcomplex}

 In this section we consider the subcomplex $P\overline{\G}$ of $P\G$ spanned by $\O$-graphs with no vertices colored by
$1_{\O}$. This subcomplex is important for calculations, as well as being neccesary to  draw the connections to $Out(F_r)$ and moduli space later on. 
  In this section we prove

\begin{proposition}\lbl{sum} The homology of $P\G$ splits as a direct sum of vector spaces
$$H_*(P\G)\cong H_*(\sy(\infty)) \oplus H_*(P\overline{\G})$$
\end{proposition}

 We first consider the subcomplex
$\mathcal B$ of
$P\G$
spanned
by graphs which \emph{only} have vertices colored by $1_\O$ (and is hence  spanned by $k$-gons).
Then $P\G   \cong   \mathcal B\oplus P\G/\mathcal B$, and Proposition~\ref{sum} will follow by showing  that $H_*(P\G/\mathcal B) \cong
H_*(P\overline{\G})$ (Lemma~\ref{modB})
and that $H_*(\mathcal B)\cong H_*(\sy(\infty))$ (Lemma~\ref{B}).

The proof of Lemma~\ref{modB} depends on understanding what happens when trivial bivalent vertices (i.e.\ vertices colored with $1_{\O}$) are added
to an edge of an $\O$ graph. For $k\geq 0$, let
$E_k$ be the oriented linear graph with $k$ internal bivalent vertices. 
Define a chain complex $\mathcal E_*$ with one generator in each dimension, corresponding to $E_k$, and boundary operator $\bdry_E$ which sums
over all edge collapses.

\begin{lemma}\lbl{edgechains} $H_0(\mathcal E_*)=\R$, and $H_i(\mathcal E_i)=0$ for $i>0$.
\end{lemma}
\begin{proof}
Collapsing any edge of $E_k$ results in a copy of $E_{k-1}$.  If we orient $E_k$ by numbering its edges (see Proposition~\ref{tree}), 
then collapsing the $j$th edge gives $(-1)^jE_{k-1}$.  Thus the boundary operator $\bdry\colon \mathcal E_j\to \mathcal E_{j-1}$ is $0$ for $j$ odd
and an isomorphism for $j$ even, giving the result.
\end{proof}

\begin{lemma}\lbl{modB}
$H_*(P\G/\mathcal B) \cong H_*(P\overline{\G})$. 
\end{lemma}
\begin{proof}
The quotient complex $\mathcal T:= P\G/\mathcal B$ is spanned by $\O$-graphs with at least one vertex
which is not colored by $1_\O$.  We filter $\mathcal T$ by the number of
vertices not colored by $1_\O$. In the associated spectral sequence
we have $$E^0_{p,q}=\F_p\mathcal T_{p+q}/\F_{p-1}\mathcal T_{p+q}.$$ 
That is $E^0_{p,q}$ is  spanned by graphs with $q$ vertices colored by $1_\O$ and $p$ other
vertices. Hence along the $p$ axis of the spectral sequence we have $E^0_{p,0} = P\overline{\G}_p$.
  The differential $d_0:E^0_{p,q}\to E^0_{p,q-1}$ is induced by the boundary operator, and is
therefore given by the map which only contracts edges incident to $1_\O$-colored vertices.
Passing to the $E^1$ term, we claim that only the $E_{p,0}$ terms persist, and that the maps
$d_1:E^1_{p,0}\to E^1_{p-1,0}$ are just the graph boundary operator, which will conclude the proof.

We break up each vertical complex 
$$\cdots E^0_{p,q}\to E^0_{p,q-1}\to E^0_{p,q-2}\to\cdots E^0_{p,0}$$
into a direct sum of chain complexes $E_p^X$, where the sum is over isomorphism classes of graphs with $p$
vertices.  An $\O$-graph is in $E_p^X$ if removing the (bivalent) $1_\O$-colored vertices from
the underlying graph results in $X$. 
Let $\Xi$ be the graph $X$ with all vertices and edges distinguished, so that $Aut(\Xi)=\{id\}$.
Then $E^X_q=E_q^{\Xi}/Aut(X)$.
 Since we are in characteristic zero, we can ignore the $Aut(X)$ action, i.e.\ $H_q(E_*^X;\R)\cong H_q(E_*^\Xi;\R)$.
Thus we must show $H_q(E_*^\Xi)$ vanishes for $q>0$. 

Since the boundary operator only affects vertices colored by
$1_\O$,
$E_*^\Xi$ breaks up as a direct sum of chain complexes $E_*^\Xi(S_1,S_2,\ldots, S_p)$, where the $S_i$ are
choices of $\O$-spider for each vertex of $\Xi$, ranging over a basis for the quotient of
vector spaces
$\Sp/<1_\O>$. Each of these chain complexes is isomorphic to $\bigotimes_{e\in E(\Xi)} \mathcal E_*$, where
$\mathcal E_*$ is the chain complex of Lemma~\ref{edgechains}.  By this Lemma and the K\"unneth formula the
homology of
$\bigotimes_{e\in E(\Xi)} \mathcal E_*$ is zero except in degree $0$.
 \end{proof}

\begin{lemma}\lbl{B}
$H_*(\mathcal B) = Prim(H_*(\sy (2\infty)))$
\end{lemma}
\begin{proof}
Direct computation shows that $H_k(\mathcal B)$ is nonzero iff $k\equiv 3 \mod 4$, in which case it is one
dimensional. This is well known to be the primitive homology of $\sy(2\infty)$, but in fact this follows if
we consider the operad which is spanned by the identity in degree 1 and is zero in higher degrees. 
Then $\A$ consists of symplecto-spiders which have $2$ legs and whose internal vertex is labeled by
$1_\O$. This is isomorphic to $\sy(2n)$ by Proposition~\ref{symplectic}. Thus Theorem~\ref{hopfiso}
implies that the primitive homology of $\a\O_\infty=\sy(2\infty)$ is given by connected $1_\O$-colored
graphs. This is precisely the complex
$\mathcal B$. \end{proof}

\subsubsection{Brackets, cobrackets, compatibility}\lbl{sec:bracket}

As we remarked earlier, the map $\phi_n$ does not establish an isomorphism of homologies.  It is not
even a chain map. However, $\G$ has a different boundary operator for which $\phi_n$ \emph{is} a chain
map.  The difference between these two boundary operators was studied in \cite{CoV} in the commutative
case, and in \cite{C} for general  cyclic operads.

\begin{definition}
Let $x$ and $y$ be half-edges of a graph $X$,   and let $\pi_{xy}$ denote the pairing of the
half-edges  $H(X)$ which differs from the standard pairing only at
$x,y,\bar x,\bar y$, where it pairs $\{x,y\}$ and $\{\bar x,\bar y\}$. If $\X$ is an $\O$-graph with underlying graph $X$, the \emph{H-boundary}
of $\X$ is 
$$\bdry_H(\X)=\sum_{x,y,x\neq \bar y} (\X^{\pi_{xy}})_{x\cup y}.$$
\end{definition}

The following proposition shows that the Lie boundary map $\bdry_n$ corresponds in fact to $2n\bdry_E + \bdry_H$.
\begin{proposition}
The following diagram commutes.
\begin{equation*}
\begin{CD} 
\O\g@>\phi_n>>\wedge\A  \\ 
@V{2n\bdry_E+\bdry_H}VV @VV{\bdry_n}V \\ 
\O\g@>\phi_n>>\wedge\A
\end{CD}
\end{equation*}
\end{proposition}

The notation $\bdry_H$ is justified by the following Lemma, which shows that $\bdry_H$ is a boundary operator in its own right.
\begin{lemma}\lbl{compat} \begin{itemize}
\item[a)]$\bdry_H^2=0$
\item[b)]$\bdry_H\bdry_E+\bdry_E\bdry_H=0$
\end{itemize}
\end{lemma}
\begin{proof}
$\phi_n$ is injective for large $n$, and $\bdry_n^2=0$, implying
$(2n\bdry_E+\bdry_H)^2=(2n)^2\bdry_E^2+2n(\bdry_H\bdry_E+\bdry_E\bdry_H)+\bdry_H^2=0$
for infinitely many values of $n$. Hence each coefficient of this polynomial in $n$ must be zero.
\end{proof}

   With respect to the disjoint union product on $\G$,
$\bdry_E$ is a derivation, but $\bdry_H$ is not.  The deviation from being a derivation defines a bracket:

\begin{definition} For $\X, {\bf Y}\in \G$, the \emph{$\O$-bracket} $[\X, {\bf Y}]$ is given by
$$[\X,{\bf Y}]=\bdry_H(\X\cdot {\bf Y})-\bdry_H(\X)\cdot {\bf Y} + (-1)^{|\X|}\X\cdot
\bdry_H({\bf Y})=\sum_{x\in \X, y\in {\bf Y}} (\X\cdot {\bf Y})_{xy},$$
where $|\X|$ is the number of vertices of $\X$.
\end{definition}

This bracket is a straightforward generalization of the bracket we defined in \cite{CoV} for the commutative
case, and coincides with the operation $\phi_2$ of \cite{C}. In \cite{C}, the following proposition is
proven.

\begin{proposition} The $\O$-bracket is graded symmetric and satisfies the graded Jacobi identities.
\end{proposition}

To obtain a graded {\it anti}-symmetric Lie bracket, we simply shift the grading (see \cite{CoV}).

The proof of the following proposition can be found in \cite{C}.
\begin{proposition}\lbl{vanish}
The $\O$-bracket vanishes canonically on $\O$-graph homology.
\end{proposition}

\begin{corollary}
$\O$-graph homology $H_*(\G)$, equipped with the differential $\bdry_H$, is a differential graded algebra (DGA).
\end{corollary}
\begin{proof}
By Lemma~\ref{compat}, $\bdry_H$ descends to the homology level. By Proposition~\ref{vanish}, it is a
derivation on the homology level. 
\end{proof}

The boundary operator $\bdry_E$ is also a co-derivation, but
the boundary operator $\bdry_H$ is not. Hence we can define
$$\theta(\X) = \Delta\bdry_H(\X)-\bdry_H\Delta(\X),$$ where $\bdry_H$ is extended to $\G\otimes\G$ as a graded
derivation. The map $\theta$ was considered in the commutative case in \cite{CoV}, and coincides with the
map $\theta_2$ of \cite{C} when the graph is connected.

Let us describe $\theta$ in more detail. A pair $x,y$ of half-edges of a graph $X$ is called
\emph{separating} if $x$ and $\bar x$ (equivalently $y$ and $\bar y$) lie in different components of
$X^{\pi_{x,y}}$. (Recall
$\pi_{xy}$ is the pairing which pairs $x$ with $y$ and  $\bar x$ with $\bar y$ and is the same as the
standard pairing everywhere else.)
Given an $\O$-graph $\X$ based on $X$, and a separating pair $x,y$ of half-edges  in $X$, write 
$\X^{\pi_{xy}} $ as a product of connected $\O$-graphs:
$$\X^{\pi_{xy}} =  \X^\prime\cdot \X^{\prime\prime}\cdot \X_1\cdot \X_2\ldots\cdot \X_m,$$
with  $x,y \in \X^\prime$ and $\bar x,\bar y \in\X^{\prime\prime}$.
Suppose $(I,J)$ is an ordered partition of $[m]$. Define the sign $\epsilon_1 (I,J)$ to satisfy the
equation
$$\X^\prime\cdot \X^{\prime\prime}\cdot \X_1\cdot \X_2\ldots\cdot \X_m= \epsilon_1(I,J) \X_I \cdot \X^\prime
\cdot \X_J
\cdot
 \X^{\prime\prime}$$ and the sign $\epsilon_2(I,J)$ to satisfy the equation
$$\X^\prime\cdot  \X^{\prime\prime}\cdot \X_1\cdot \X_2\ldots\cdot \X_m= \epsilon_2(I,J) \X_I \cdot 
\X^{\prime\prime}
\cdot \X_J
\cdot \X^\prime.$$ Now we define\quad $\theta (\X)$
$$=\sum_{x,y \text{ separating}}\sum_{I,J} \epsilon_1 (I,J) (\X_I \cdot \X^\prime_{x\cup y})
\otimes (\X_J
\cdot \X^{\prime\prime})
+\epsilon_2(I,J) (\X_I \cdot \X^{\prime\prime})\otimes (\X_J \cdot \X^\prime_{x\cup y}).$$
  It is straightforward to verify that this agrees with the previous definition of $\theta$.

The following proposition is proven in \cite{C}.

\begin{proposition}  The map $\theta\colon\G\to\G\otimes\G$ restricts to a map $P\G\to P\G\otimes P\G$, 
and is a (graded) cobracket on $P\G$. 
\end{proposition}

Conjecturally, $\theta$ is a cobracket on the entire complex. In \cite{CoV}, we proved this for the commutative
case, and the proof extends to the associative case. However, it does not extend to the Lie case.

The following was noted without proof in \cite{C} in the connected case.
\begin{proposition}\lbl{vanish2}
The map $\theta$ vanishes canonically on the homology level. 
\end{proposition}
\begin{proof}
Define a map $T\colon \G\to\G\otimes \G$ as
$$T(\X)=\sum_{x,y \text{ separating}}\sum_{I,J} \epsilon_1 (I,J) (\X_I \cdot \X^\prime) \otimes
(\X_J
\cdot  \X^{\prime\prime}),$$
where $\X$ is written as a product of connected graphs as above. 
For notational ease, write $T(\X)=\sum_{x,y\text{ separating}} T_{xy}(\X)$.
We claim that $\theta = \bdry_E T - T
\bdry_E$, i.e.\ $\bdry_E T= T
\bdry_E+ \theta$.  To see this, consider the term $\bdry_E T_{xy}(\X)$ of $\bdry_E T(\X)$, which contracts all edges of 
$T_{xy}(\X)$, one at a time. If the edge $e$ is not $x\cup y$ or $\bar x\cup \bar y$, then we have
$(T_{xy}(\X))_e=T_{xy}(\X_e)$. In this way we pick up all terms of $T \bdry_E(\X)$. If the edge is 
$x\cup y$ or $\bar x\cup\bar y$, then we pick up the two terms in $\theta (\X)$.

Now suppose that $\X$ is a cycle. Then the above claim implies that $\theta(\X) = \bdry_E T(\X)$, 
and is hence trivial on the homology level.
\end{proof}

\begin{corollary}
$\O$-graph homology $H_*(\G)$ is a differential graded Hopf algebra, with differential $\bdry_H$.
\end{corollary}

\begin{corollary}
The map  $\bdry_H$ is a graded differential on $H_*(P\G)$ and $H_*(P\overline{\G})$.
\end{corollary}
\begin{proof} 
Since $\bdry_H$ is a derivation on the homology level, it takes primitives to primitives, and so induces a differential
on $Prim H_*(\G)=H_*(P\G)$. 
This is \emph{not}
true on the chain level, as there exist connected graphs with disconnected terms of $\partial_H$. 

Also any graph with a vertex colored by $1_\O$ will get mapped by $\partial_H$ to a sum of graphs each with a
vertex colored by $1_\O$ implying that $\partial_H$ is a differential on $H_*(P\overline{\G})$.
If the pair of half-edges being contracted in $\partial_H$ is not adjacent to a $1_\O$-colored vertex,
then the $1_\O$ vertex survives as claimed. Otherwise, consider a $1_\O$ colored vertex. It has 
two adjacent half-edges $h_1,h_2$. Given another half-edge $k$, we claim that the
$h_1,k$ term cancels with the $h_2,\bar k$ term. In other words all the summands coming from contracting
edges adjacent to $1_\O$ colored vertices cancel! Clearly these two terms give isomorphic graphs,
and checking the orientation one sees that they have opposite sign.
\end{proof}

The following proposition is proven in \cite{CoV} for the commutative case, and as remarked in \cite{C}, the
proof carries over to the general cyclic operad case.
It explains our interest in graphs with no separating edges.

\begin{proposition}\lbl{bialgebra}  Bracket and cobracket form a compatible Lie bi-algebra structure on the subcomplex of
connected, 1-particle irreducible graphs (graphs with no separating edges).
\end{proposition}

\section{The Lie case}

\subsection{The forested graph complex and the Lie graph complex}

We now specialize to the Lie operad. 
 Recall that $\O[m]$ is the vector
space spanned by all
 rooted vertex-oriented binary trees with $m$ numbered
leaves, modulo the subspace
spanned by all anti-symmetry and IHX relators.
A \emph{vertex
orientation} is
a choice of cyclic order at each trivalent vertex of a binary tree, and is equivalent to specifying a
planar embedding up to isotopy.
 The composition rule in the operad attaches the root of
one tree to a leaf of another, eliminates the resulting bivalent vertex, and then suitably renumbers
the remaining leaves. 

The space of $\O$-spiders, $\Sp$, is spanned by vertex-oriented binary trees with $m$ numbered leaves,
but with no particular leaf designated as the root. Mating is defined by gluing two such trees together
at a leaf, eliminating the resulting bivalent vertex, then renumbering the leaves suitably. 

Thus a basic $\O$-graph for the Lie operad is a fairly complicated object: an oriented graph with a
vertex-oriented trivalent tree attached at each vertex, modulo AS and IHX relations.  However, the picture can be
considerably simplified  since the orientations of the trees at the vertices cancel to a large degree with
the orientation of the graph itself. We now describe this simplification.

Let $X$ be a finite graph.  A \emph{ forest} in $X$ is a sub-graph which contains no cycles.  The 
connected components of a forest are trees, where we allow a tree to consist of a single vertex.  An
\emph{ orientation} of a forest is  given by an ordering of its edges (if any); interchanging the 
order of any two edges  reverses
the orientation of the forest.

\begin{definition} A \emph{forested graph} is a pair $(G,\Phi)$, where $G$ is a finite  
connected trivalent  graph  and $\Phi$ is an oriented forest
  which contains all vertices of $G$.
\end{definition}
An example of a forested graph with three trees is shown in Figure~\ref{fgraph}.
\begin{figure}[ht!]
\begin{center}
\includegraphics[height=2cm]{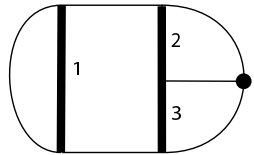} 
\end{center}
\caption{Forested graph}\lbl{fgraph} 
\end{figure}

We denote by $\V_k$  the vector space spanned by all forested graphs whose forest contains exactly $k$
trees, modulo the relations
$(G,\Phi)  = -(G,-\Phi)$.   

Let $(G,\Phi)$ be a forested graph, and let $e$ be an edge of $\Phi$.  Collapsing $e$ produces a new pair
$(G_e,\Phi_e)$, where  $G_e$ has one 4-valent vertex.  There are exactly two 
other forested graphs (which may be isomorphic!), 
$(G^{'},\Phi{'})$ and $(G^{''},\Phi{''})$, and edges $e'\in \Phi',e{''}\in \Phi{''}$  which produce the
same pair $(G_e,\Phi_e)$ (see Figure~\ref{fIHX}). 
\begin{figure}[ht!]
\begin{center}
\includegraphics[height=6cm]{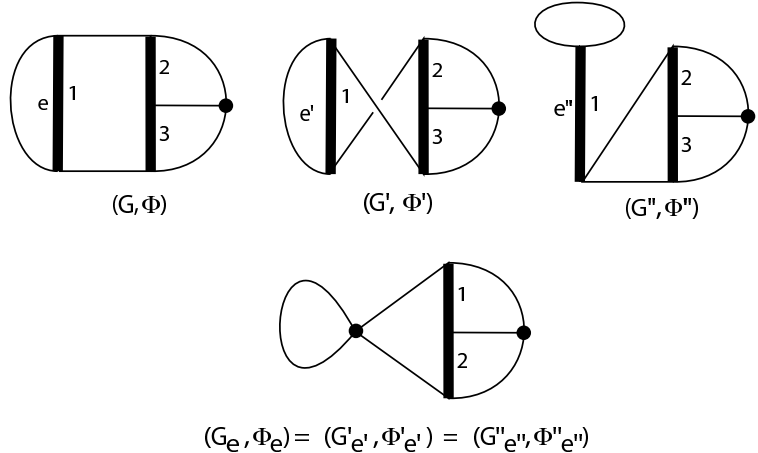}
\end{center}
\caption{Blow-ups of $(G_e,\Phi_e)$}\lbl{fIHX}
\end{figure}

\noindent The vector    
          $$(G,\Phi)+(G^{'},\Phi{'})+(G^{''},\Phi{''})$$
is called the  \emph{ basic IHX relator} associated to $(G,\Phi,e)$.  

Let $IHX_k$ be the subspace of $\V_k$ spanned by all basic IHX relators, and
define
$\C_k$ to be the quotient space $\V_k/IHX_k$.

\begin{proposition}\lbl{Ck} $\C_k$ is
naturally isomorphic as a vector space to 
$P\overline\G_k$, where $\O$ is the Lie operad.
\end{proposition}

\begin{proof}
We define a map $\V_k\to P\overline\G_k$ as follows.  Starting
with a forested graph $(G,\Phi)$, collapse each component of $\Phi$ to a point to produce a graph $X$. 
For each vertex $v$ of $X$, let $T_v$ be an
$\epsilon$-neighborhood of the preimage of $v$ in $G$, so that $T_v$ is a binary tree whose leaves
 are identified with
the half-edges $H(v)$ adjacent to $v$  (see Figure
~\ref{correspondence}).  To specify an $\O$-graph
$\X=(X,\{T_v\})$, it remains only to determine the orientations on $X$ and on the trees $T_v$.  
\begin{figure}[ht!]
\begin{center}
\includegraphics[width=.9\hsize]{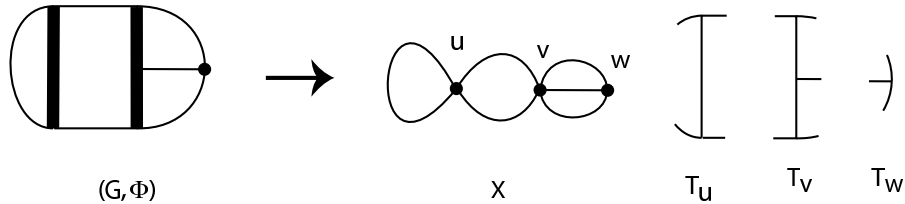}
\end{center}
\caption{Correspondence between forested graphs and Lie graphs}\lbl{correspondence}
\end{figure}
  In the language of section~\ref{sec:orientation}, this information can be described as a unit vector in 
$$\det(\R V(X))\otimes \bigotimes_{e\in E(X)} \det \R H(e)\otimes \bigotimes_{v\in
V(X)}(\otimes_{u\in V(T_v)}
\det \R H(u))$$
That is, it is specified by a numbering of the vertices of $X$, an orientation on each edge of $X$,
and a vertex-orientation on each tree $T_v$.

Applying Propositions~\ref{orient} and \ref{orient2} to the first two tensor factors and 
applying Corollary ~\ref{tree} to the third tensor factor, this is canonically isomorphic to 
\begin{align}
\det(\bigoplus_{|H(v)| even} \R v)\otimes\bigotimes_{v\in V(X)}\det \R H(v)\otimes
  \bigotimes_{v\in
V(X)}\det (\R E(T_v)) \lbl{eqtn}
\end{align}
In other words, an orientation can be  specified by an ordering of the even valence vertices of $X$, an ordering of the
half-edges incident to each vertex of $X$, and an ordering of the edges of each tree.

The forest $\Phi$ is the union of the \emph{internal} edges of all of the $T_v$.
Denote the internal edges of $T_v$ by $\Phi_v$, so that $\Phi=\cup_v \Phi_v$. 
Since $H(v)$ is identified with the leaves of $T_v$, the Partition lemma implies 
$\det(\R E(T_v))\cong\det(\R H(v))\otimes\det(\R E(\Phi_v))  
 $.
(Actually, the Partition lemma tells you to also order the set of odd subsets of $E(T_v)$, but since 
the number of leaves and the number of interior edges have opposite parity,
this would be ordering at most one object and that is not extra information.)
Hence,  
  $\det(\R H(v))\otimes \det(\R E(T_v)) \cong
    \det(\R E(\Phi_v))$, so that expression~ (\ref{eqtn}) becomes
\begin{align*}
&  \det(\bigoplus_{|H(v)| even} \R v)\otimes\bigotimes_{v\in V(X)}\det (\R E(\Phi_v))\\
& = \det(\bigoplus_{|E(\Phi_v)| odd} \R v)\otimes\bigotimes_{v\in V(X)}\det (\R E(\Phi_v))\\
& \cong \det(\bigoplus_{e\in \Phi} \R e),
\end{align*}
where the last isomorphism is again given by the Partition lemma. 
Thus the orientation data on $\X=(X,\{T_v\})$ is equivalent to an ordering of the edges of $\Phi$. 

This defines the map $\V_k\to P\overline\G_k$.  We claim that a basic IHX relator
$(G,\Phi)+(G',\Phi')+(G'',\Phi'')$ maps to
$(X,\{T_v\})-(X,\{T'_v\})+(X,\{T''_v\})= (X,\{T_v-T'_v+T''_v\})=0,$ so the map factors through $\C_k$.
This is obvious except for the signs of the terms.  To check these, we need to carefully  consider 
 the equivalence of orientations described above.
 The essential step is the translation from vertex
orienting a tree to ordering the tree's edges. 
We illustrate the argument by doing the case of a tree with $5$ edges, and  leave the general case to the reader.
Suppose we have
\begin{center}
\includegraphics[height=1.5cm]{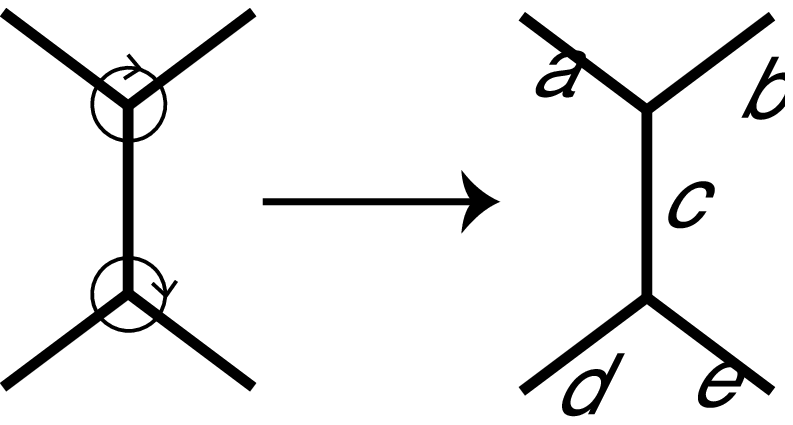} 
\end{center}
where $\{a,b,c,d,e\}=\{1,2,3,4,5\}$ (not neccesarily in that order).
Then, depending on where we place the tree in the plane,  we have 
\begin{center}
\includegraphics[height=3cm]{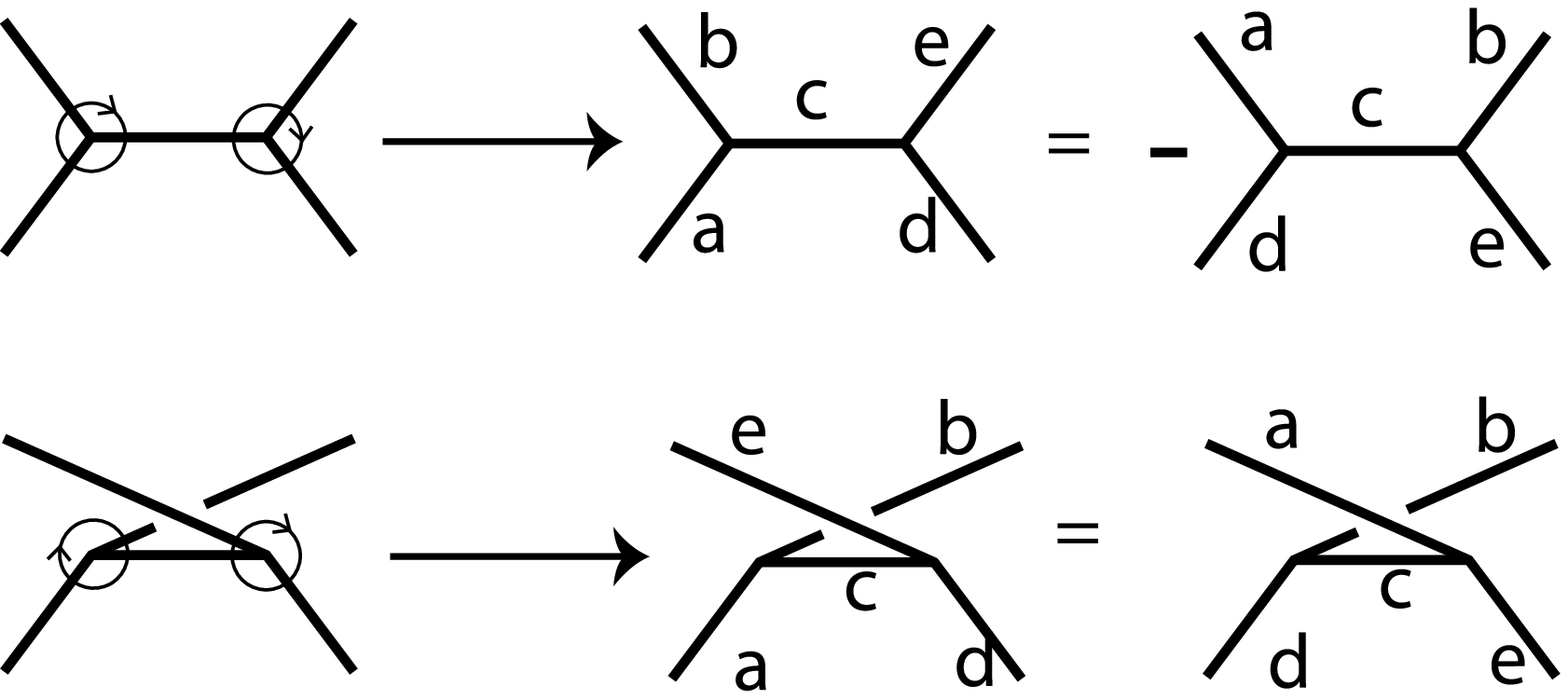}
\end{center}
Thus the $I-H+X$ relator of vertex oriented trees gets mapped to an   $I+H+X$ relator of trees with
ordered edges. 

To see that the map $\C_k\to P\overline\G_k$ is an isomorphism, define an inverse $P\overline\G_k\to \V_k\to \C_k$ as follows:
If $(X,\{T_v\})$ is a basic $\O$-graph,
produce a
trivalent graph $G$ from this data by replacing each vertex $v$ of $X$ by 
$T_v$, using the identification of the leaves of $T_v$ with $H(v)$ as gluing instructions  (see Figure
~\ref{correspondence}). The union of the interior edges of the $T_v$ forms a forest
$\Phi$ in
$G$, and we have seen that the orientation data on   $(X,\{T_v\})$ determines an orientation of $\Phi$.
\end{proof}

We have identified the spaces $\C_k$ with the chain groups $P\overline\G_k$. Under this
identification, it can be checked that the boundary map
$\bdry_E$ is induced by the map on $\V_k$ given by
$$\bE(G,\Phi)=\sum (G,\Phi\cup e),$$
where the sum is over all edges $e$ of $G-\Phi$ such that $\Phi\cup e$ is still a forest.
(Note that this only happens if the initial and terminal vertices of $e$ lie in different
components of
 $\Phi$, so that $\Phi\cup e$ has $k-1$ components.) The orientation of $\Phi\cup e$ is
determined by ordering the edges of $\Phi$ with the labels $1,\ldots, k$ consistent with its
orientation, and then labeling the new edge $e$ with $k+1$.

\begin{remark} As defined above, $\bdry_E$ is a boundary operator on the chain complex $\V_*$ of forested graphs.  There is also a coboundary map
$\delta\colon \V_k\to \V_{k+1}$, defined by the formula
$$\delta(G,\Phi)=\sum_{e_i\in \Phi} (-1)^i(G, \Phi-e_i),$$
where the orientation on $\Phi-e_i$ is induced from the ordering of the edges of $\Phi$.  However, this
coboundary does not preserve the IHX-subspaces, so does not induce a coboundary operator on $\C_*$.  Of
course, we can get a coboundary operator by dualizing: $\delta=\bE^*\colon \C_k^*\to \C_{k+1}^*$, but
there is no natural basis of $\C_k$, and hence no natural identification of $\C_k$ with
$\C_k^*$.
\end{remark}

\subsubsection{Bracket and cobracket in the forested graph complex}
In this section we describe how the bracket and cobracket of section~\ref{sec:bracket} translate
to forested graphs.

Let $(G_1,\Phi_1)$ and $(G_2,\Phi_2)$ be forested graphs,  let  
$x$ be a half-edge in
$G_1-\Phi_1$  and let  
$y$ be a half-edge in
$G_2-\Phi_2$.  We  form a new trivalent graph $(G_1\cdot G_2)^{\pi_{xy}}$, and forest it 
by the image of $\Phi_1\cup \Phi_2$  plus the new edge  $xy$.  
The orientation  is given by shifting the numbering the edges of $\Phi_2$ to lie after $\Phi_1$,
and letting $x\cup y$ come after that.
 The bracket is then defined by summing over  all possible pairs
$x$ and
$y$:
$$[(G_1,\Phi_1), (G_2,\Phi_2)]=\sum_{x,y}((G_1\cdot G_2)^{\pi_{xy}},\Phi_1\cup\Phi_2\cup x\cup y)$$

The cobracket is defined as follows. Let $(G,\Phi)$ be a forested graph, and let $x,y$ be two
half-edges not in $\Phi$   and not adjacent to the same tree of $\Phi$.   such that $G^{\pi_{xy}}$ has
two connected components $G_1$ and $G_2$, where
$x\cup y$ is in $G_1$. The pair of half-edges $\{x,y\}$ is said to be \emph{separating}.
Let $\Phi_i=\Phi\cap G_i$. Add $x\cup y$ to the forest $\Phi_1\cup\Phi_2$,
ordered after everything else.
 Adjust the representative of the orientation so that
the numbering of $\Phi_1\cup x\cup y$ precedes that of $\Phi_2$. Form the symmetric product
$(G_1,\Phi_1\cup x\cup y)\odot (G_2,\Phi_2)$, where the numbering of $\Phi_2$ is shifted down
by the number of edges in $\Phi_1\cup x\cup y$. The cobracket is given by summing over all
such pairs $\{x,y\}$:
$$\theta (G,\Phi) = \sum_{\{x,y\}\text{ separating}} (G_1,\Phi_1\cup x\cup y)\odot (G_2,\Phi_2)$$

\subsection{A filtration of Outer Space, and the associated complex}

The group $Out(F_r)$ of outer automorphisms of a free group of rank $r$ acts cocompactly on a contractible simplicial
complex $K_r$, known as the spine of Outer space.  For an introduction to Outer space and its spine, see \cite{V2}.   Point stabilizers of this action are
finite, so that the quotient $Q_r=K_r/Out(F_r)$ 
 has the same rational cohomology as $Out(F_r)$. In this section we define a filtration on $K_r$ and use the spectral
sequence of this filtration to prove the following theorem:

\begin{theorem}\lbl{homology}
Let $\C^{(r)}$ denote the subcomplex of $\C$ spanned by (connected) forested graphs of rank $r$.
Then $H_k(\C^{(r)})\cong H^{2r-2-k}(Out(F_r))$. 
\end{theorem} 

 Recall that vertices of $K_r$ are
pairs
$(g,X)$, where
$X$ is a connected graph with all vertices at least trivalent, and $g$ is a homotopy equivalence from a 
standard  rose (wedge of circles)  $R_r$ to $X$; thus $g$ gives an
identification of $\pi_1(X)$ with the  group $F_r= \pi_1(R_r)$.  The $k$-simplices of $K_r$, for $k\geq 1$, can be identified with
chains $$\emptyset=\Phi_0\subset\Phi_1\subset\ldots\subset\Phi_k$$  of forests in a marked graph $(g,X)$; the $i$th vertex of the $k$-simplex is
obtained from
$(g,X)$ by collapsing each edge of
$\Phi_i$ to a point.   $Out(F_r)$ acts by changing the marking $g$, and the stabilizer of
$(g,X)$ is isomorphic to the group of automorphisms of the graph $X$.

We define a filtration on $K_r$ as follows:
\begin{itemize}
\item {} $\F_0K_r$ consists of  all vertices $(g, X)$ with $X$ a (connected) trivalent graph.  Note that
a trivalent graph $X$ with fundamental group isomorphic to $F_r$ has
$2r-2$ vertices.
\item{}$\F_iK_r$ is the sucomplex of $K_r$ spanned by $\F_{i-1}K_r$ together with all vertices $(g,X)$
such that
$X$ has $2r-2-i$ vertices.  
\end{itemize}
We have  $\F_0K_r\subset \F_1K_r\subset \ldots \subset
\F_{2r-3}K_r=K_r$, where the $i$th subcomplex $\F_iK_r$ has dimension $i$.  
Since the action of 
$Out(F_r)$ simply changes the markings $g$, it preserves the filtration, and so induces a filtration 
$\F_0Q_r\subset \F_1Q_r\subset \ldots \subset
\F_{2r-3}Q_r=Q_r$ on the quotient $Q_r$.  The cohomology spectral sequence of this filtration has
$$E_1^{p,q}=H^{p+q}( \F_pQ_r, \F_{p-1}Q_r)$$
and converges to the cohomology of $Q_r$.  In particular, if we take trivial real coefficients, this
spectral sequence converges to the rational cohomology of $Out(F_r)$.  All cohomology groups will be assumed
to have trivial rational coefficients unless otherwise specified.
 
\begin{proposition}\lbl{zero} $ H^{p+q}( \F_pQ_r, \F_{p-1}Q_r)=0$ for $q\neq 0$
\end{proposition}

\begin{proof} Marked graphs in $\F_pK_r$ which are not in $\F_{p-1}K_r$ are those where $\sum_v (|v|-3)$ is exactly $p$, where
$|v|$ is the valence of the vertex $v$.  The link in $\F_{p-1}K_r$ of
such a marked graph $(g,X)$ is the join of the blow-up complexes of the vertices of $X$ of valence $>3$.  

A simplicial complex is called \emph{ $i$-spherical} if it is $i$-dimensional and homotopy equivalent to a wedge of $i$-spheres.
The blow-up complex of
a vertex of valence
$k$ is $(k-4)$-spherical (see, e.g. \cite{V1}, Theorem 2.4).  Therefore the link of $(g,X)$ in $\F_{p-1}K_r$ is
$(p-1)$-spherical.  Since graphs in $\F_pK_r-\F_{p-1}K_r$ are not connected by edges, the entire complex
$\F_pK_r/\F_{p-1}K_r$ is $p$-spherical. Therefore   $ H^{p+q}(\F_pK_r,\F_{p-1}K_r)=0$ for $q\neq 0$. 

Since $Out(F_r)$ preserves the filtration, the homology of $\F_pQ_r/ \F_{p-1}Q_r$ is the homology of
$\F_pK_r/\F_{p-1}K_r$ modulo the action of $Out(F_r)$.  The stabilizer of each graph in $\F_pK_r$ is finite, so the (reduced)
homology of the quotient $\F_pQ_r/\F_{p-1}Q_r$ vanishes in dimensions other than $p$, as was to be shown.
\end{proof}

\begin{proof}[Proof of Theorem]

By Proposition~\ref{zero}, the spectral sequence of the filtration $\F_*Q_r$ of $Q_r$ collapses to the cochain complex
\begin{align*}
0\to H^0(\F_0Q_r)\to H^1(\F_1Q_r, \F_0Q_r)\to &H^2( \F_2Q_r, \F_1Q_r)\to \\
&\ldots\to H^{2r-3}(
\F_{2r-3}Q_r,\F_{2r-4}Q_r)\to 0,
\end{align*}
which computes the cohomology of $Out(F_r)$.

To identify the terms $H^p( \F_pQ_r, \F_{p-1}Q_r)$ it is helpful to recall the description of $K_r$ as a
cubical complex (see \cite{HV}).    The simplices of $K_r$ organize themselves naturally into cubes $(g,X,\Phi)$, one for each forest $\Phi$ in a
marked graph $(g,X)$ (see Figure~\ref{cube}).  An ordering on the edges of $\Phi$ determines an orientation of the cube.

\begin{figure}[ht!]\begin{center}
\includegraphics[height=7cm]{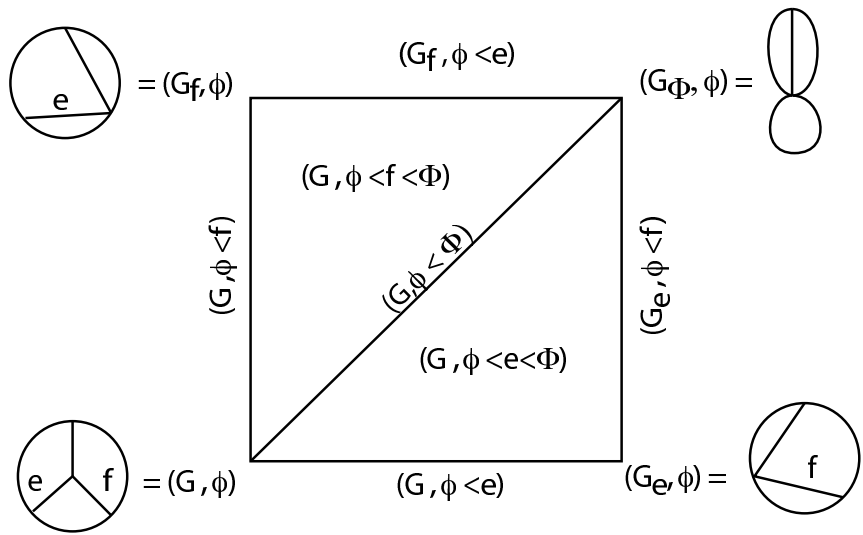}
\caption{The cube $(G,\Phi)$ for $\Phi=\{e,f\}$}\lbl{cube}
\end{center}
\end{figure}

 The stabilizer of a cube $(g,X,\Phi)$ consists of automorphisms
of $X$ which preserve $\Phi$ (see \cite{SV}).  The quotient of a cube by its stabilizer is a cone on a rational
homology sphere, or on a rational homology disk if there is an automorphism of $(X,\Phi)$ which reverses
orientation.  

Now\quad $H^p(\F_pQ_r,\F_{p-1}Q_r)$
$$=C^p(\F_pQ_r,\F_{p-1}Q_r)/im(\delta\colon C^{p-1}(\F_pQ_r,\F_{p-1}Q_r)\to
C^p(\F_pQ_r,\F_{p-1}Q_r)).$$  Note that
$\F_{p-1}Q_r$ is $(p-1)$-dimensional, and $\F_pQ_r$ has one $p$-dimensional cube quotient for every pair $(G,\Phi)$ such that $G$ is
trivalent and  $\Phi$ has $p$ edges.  Therefore
$$C^p(\F_pQ_r,\F_{p-1}Q_r)=C^p(\F_pQ_r)=( \R\{(G,\Phi)\}/AS)^*,$$
where 	$AS$ is the subspace generated by the antisymmetry relations $(G,\Phi)=-(G,-\Phi)$.

Each $p$-cube $(g,G,\Phi)$ in $\F_pK_r$ has $2^p$ codimension one faces of the form $(g,G,\Phi-e)$, which are in
$\F_{p-1}K_r$.  Opposite each is a codimension one face $(g,G,\Phi-e)^{op}=(g_e,G_e,\Phi_e)$.  
The relative cochain group $C^{p-1}(\F_pK_r,\F_{p-1}K_r)$ consists of functions on the  $(p-1)$-cubes of $\F_{p}K_r$ which
vanish on $\F_{p-1}K_r$, i.e.\
 $C^{p-1}(\F_pK_r,\F_{p-1}K_r)$ has as a basis the characteristic functions of the $(p-1)$-cubes $(g_e,G_e,\Phi_e)$ of the
second type.   Each such $(p-1)$-cube is a face of exactly three $p$-cubes, $(g,G,\Phi)$, $(g',G', \Phi')$ and
$(g{''},G{''},\Phi{''})$, with $G_e=G'_{e'}=G{''}_{e{''}}$, so that $\delta(g_e,G_e,\Phi_e)=(g,G,\Phi)+(g',G',
\Phi')+(g{''},G{''},\Phi{''})$.  Passing to the quotient modulo $Out(F_r)$, we get the IHX relator
$$\delta(G_e,\Phi_e)=(G,\Phi)+(G',\Phi')+(G{''},\Phi{''}),$$
 showing that in the quotient,
$im(\delta)$ is spanned by the (characteristic functions of) the IHX relators $IHX(G,\Phi,e)$.  

Thus\quad$H^p(\F_pQ_r,\F_{p-1}Q_r)$
$$=C^p(\F_pQ_r,\F_{p-1}Q_r)/im(\delta )= 
(\R\{(G,\Phi)\}/AS)^*/IHX^*=\C_{2r-2-p}$$
The last isomorphism follows since $\R\{(G,\Phi)\}/AS$ and $IHX$ are canonically isomorphic to their duals,
 even though the quotient
is not.
This identifies the cochain groups $E_1^{p,0}$ above with the chain groups of the forested graph complex.  The $d^1$ map in the
spectral sequence is induced by the coboundary map $$H^{p-1}(\F_{p-1}Q_r,\F_{p-2} Q_r)\to H^p(\F_pQ_r,\F_{p-1}Q_r)$$ 
which maps
$$(G,\Phi)\mapsto\sum (G,\Phi\cup e),$$
where the sum is over all edges $e$ in $G$ such that $\Phi\cup e$ is a forest;
this is the $\bdry_E$ map in the forested graph complex.
\end{proof}

\subsection{The\,subcomplex\,of\,connected\,one-particle\,irreducible\,graphs}

A graph $X$ is called \emph{ one-particle irreducible} if removing an edge does not change the number of
connected components, i.e.\ if $X$ has no  separating edges.  We have

\begin{proposition} The subcomplex $IK_r$ of $K_r$ spanned by (connected) one-particle irreducible graphs is
an equivariant deformation retract of $K_r$.
\end{proposition}

\begin{proof}   The deformation retraction acts by uniformly shrinking all separating edges to points in all
graphs representing vertices of $K_r$.  (see \cite{CuV}).
\end{proof}

Thus we may compute the rational cohomology of $Out(F_r)$ using $IK_r$ instead of $K_r$.  The filtration $\F_*K_r$ restricts to a filtration on
$IK_r$, and of the quotient ${IQ}_r$ of $IK_r$ by the action of $Out(F_r)$, giving a spectral sequence 
$$E_1^{p,q}=H^{p+q}( \F_p{IQ}_r, \F_{p-1}{IQ}_r)$$
The fact that this sequence degenerates to a chain complex computing the homology of $Out(F_r)$ is guaranteed by 
the following proposition

\begin{proposition} $H^{p+q}( \F_p{IQ}_r, \F_{p-1}{IQ}_r)=0$ for $q\neq 0$
\end{proposition}

\begin{proof} As in the proof of Proposition~\ref{zero},  what is needed is to prove that the link of a vertex in
$\F_pIK_r-\F_{p-1}IK_r$ intersects
$\F_{p-1}IK_r$ in a
$(p-1)$-spherical subcomplex. This link is called $K_{>(g,G)}$ in \cite{V1}, and Corollary 3.2 of Theorem 3.1 
of that paper proves that it is $(p-1)$-spherical.  
\end{proof}

Now let $\O$ be the Lie operad,
and consider the complex $PI\overline\G^{(r)}$ of $P\overline\G^{(r)}$ which is the quotient of $P\overline\G^{(r)}$
by graphs with separating edges, which may be internal to vertices.
Proposition ~\ref{bialgebra} implies that the bracket and cobracket 
form a compatible bi-algebra structure on this complex. (One must check that the bracket takes 
graphs with internal separating edges to graphs with internal separating edges, and check a similar statement for the cobracket.) 
 The proof of Proposition ~\ref{Ck} restricts to
give an isomorphism of the subspace $I\C_k$ of $\C_k$ spanned by 1-particle irreducible forested
graphs with $PI\overline\G_k$.  The proof of Theorem ~\ref{homology}  also restricts, and the
combination gives 

\begin{theorem}
$H_k(PI\overline\G^{(r)})\cong H_k(I\C^{(r)})\cong H^{2r-2-k}(Out(F_r))$.\qed 
\end{theorem} 

  This theorem, together with Proposition~\ref{bialgebra}, shows  there is a chain complex that
computes $H^*(Out(F_r);\R)$ which carries a highly non-trivial Lie bialgebra structure.  This bialgebra structure
vanishes at the level of homology, so that $\oplus_r H^*(Out(F_r);\R)$ is the primitive part of a
differential graded Hopf algebra.  

\section{The associative case}
\subsection{The associative graph complex and the forested ribbon graph complex}
We now make a few remarks concerning the associative case, i.e when the operad $\O$ is the {\it associative
operad}. In this case
 $\O[m]$ is $m!$-dimensional, with a basis consisting of of rooted planar trees with one interior vertex and
$m$ numbered leaves. An equivalent 
description, closer to the description of the Lie operad, is that $\O[m]$ is spanned by rooted planar \emph{binary}
trees with
$m$ numbered leaves,
modulo anti-symmetry and associativity (or ``IH") relations.

$\O$-spiders in this case 
 are planar spiders with one vertex and $m$ numbered legs.  Planarity of the spider can be thought of as a
 cyclic ordering on the legs, so that  an $\O$-graph can be characterized as an oriented graph $X$ together
with a cyclic ordering
 of the edges incident to each vertex.  Such objects are known as a ``ribbon
graphs," and have been studied in  conjunction with mapping class groups of punctured surfaces (see, e.g., \cite{Penner}).  

For $S$ an oriented compact surface with boundary, we denote by $\Gamma(S)$  the mapping class group  of isotopy classes of
homeomorphisms of
$S$  which preserve the orientation.  Homeomorphisms and isotopies need not be the identity on the boundary, and in particular
homeomorphisms may permute the boundary components.  
 
\subsubsection{Surface subcomplexes of $P\overline{\G}$}
Given a (connected) ribbon graph in $P\overline{\G}$, one can ``fatten" it to a unique compact, connected, oriented surface in such a way
that the cyclic orientation at each vertex is induced from the orientation of the surface. 
Contracting an edge will not change the oriented surface. Thus the connected, reduced associative graph
complex breaks up into a direct sum over surfaces $S$:
$$P\overline{\G} = \bigoplus_S \left(P\overline{\G}\right)_S.$$
One can think of a ribbon graph as an equivalence class of embeddings of a graph into a surface $S$ which induce an 
isomorphism on $\pi_1$; here two
embeddings are considered equivalent if they differ by   a homeomorphism of the surface. 

\subsubsection{Forested ribbon graphs}
A \emph{forested ribbon graph}, $(G,\Phi)$
 is a connected trivalent ribbon graph $G$ and an oriented forest $\Phi\subset G$ that contains all the vertices.
Let $\W$ denote the vector space spanned by forested ribbon graphs, modulo the anti-symmetry relation $(G,\Phi)=-(G,-\Phi)$
. Let $\W_S$ denote the subcomplex of $\W$ spanned by 
forested graphs $(G,\Phi)$ where $G$ thickens to the surface $S$. 

Recall that
an IHX relation comes from blowing up a $4$-valent vertex in all possible ways. Similarly, in the case of a ribbon
graph, we can define an ``IH" relation which is the sum of blowing up a $4$-valent vertex in the two possible ways consistent with the 
ribbon structure. Now define $\fr$ to be the quotient of $\W$ by IH relations, and $\fr_S$ to be the subcomplex which thickens
to $S$. 

\begin{proposition}\lbl{forested ribbon}
$\fr_S$ is naturally isomorphic to $\left(P\overline{\G}\right)_S$.
\end{proposition}
\begin{proof}
We use the description of the associative operad as a quotient of binary trees by IH relations.
Define a map $\fr_S\to \left(P\overline{\G}\right)_S$ as in the Lie case by collapsing the forest, and coloring
the resulting vertices by $\epsilon$-neighborhoods of preimages.
As in the Lie case, the orientation of the forest translates to a vertex-orientation of the tree at each vertex plus an orientation 
in the usual sense of the underlying graph. Each tree is already canonically vertex-oriented by the ribbon structure, so that
we can compare the induced vertex-orientation with the canonical one, incurring a plus or minus sign.

One must check that $IH$ relations get sent to $IH$ relations and that this translation is a chain map, the
only issue in both cases being the sign.
\end{proof}

\subsection{Ribbon subcomplexes of Outer Space}

Fix an oriented surface $S$ with non-empty boundary, and consider the set of all isotopy classes of embeddings of ribbon graphs $(X,\{or_v\})\to S$,
where the embedding respects the cyclic orientations at vertices, and induces an isomorphism on the fundamental group. The set of all such isotopy
classes of embeddings forms a simplicial complex
$L_S$, where an edge corresponds to collapsing a forest in $X$.  In fact, if we forget the orientations $or_v$, this is just a subcomplex of the
spine $K_r$ of Outer space.  To see, this, we identify $\pi_1(S)$ with the fundamental group of a rose $R_r$
which is a deformation retract of $S$; then the marking
$g$ on $X$ is a homotopy inverse to the embedding followed by the retraction.

This subcomplex of $K_r$ is contractible.  This can be seen directly from the proof that Outer space is
contractible (for details, see
\cite{Horak}), or via the identification of this complex with a deformation retract of the Teichm\"uller space of the punctured surface
(\cite{Penner}).  The mapping class group  $\Gamma(S)$ of $S$ is naturally a subgroup of $Out(F_r)$, namely,
$\Gamma(S)$ is the stabilizer, in
$Out(F_r)$, of the set of cyclic words represented by small loops around the punctures (see \cite{ZVC}).  The
action of $Out(F_r)$ on $K_r$ restricts to an action of
$\Gamma(S)$ on
$L_S$, so the quotient of
$L_S$ by this action, denoted $Q_S$ is a rational
$K(\Gamma(S),1)$.   The filtration of $K_r$ by the number of vertices in the graph restricts to a filtration
$\F_pL_S$ of $L_S$ which is invariant under the action of
 $\Gamma(S)$, so the cohomology spectral sequence of the quotient filtration converges to the cohomology of $\Gamma(S)$. We have
$$E_1^{p,q}=H^{p+q}( \F_pQ_S, \F_{p-1}Q_S)$$ 

As before, we need the analog of Proposition~\ref{zero}:
\begin{proposition} $H^{p+q}( \F_pQ_S, \F_{p-1}Q_S)=0$ for $q\neq 0$
\end{proposition}

\begin{proof}  In this case, the link of a vertex in $\F_pL_S-
\F_{p-1}L_S$ intersects $\F_{p-1}L_S$ in a single
sphere, of dimension $(p-1)$. This is basically a consequence of the fact that the ``IHX" relation is just an
``IH" relation when the graph is restricted to lie on a surface, i.e.\ when it is necessary to preserve the
cyclic orderings at vertices.
\end{proof}
  
\begin{theorem}
Let $S$ be a surface with fundamental group $F_r$.
Then $$H^{2r-2-k}(\Gamma(S);\R)\cong H_k(\fr_S) \cong H_k((P\overline{\G})_S)$$
\end{theorem}
\proof
By the above remarks $H_*(\Gamma(S);\R)$ is computed by the chain complex
$$\ldots\to H^p(\F_pQ_S, \F_{p-1}Q_S)\to H^{p+1}( \F_{p+1}Q_S ,\F_pQ_S )\to\ldots$$
As in the case of $Out(F_r)$, this can be described as  a quotient of the vector space spanned by forested ribbon graphs
$$(\R\{(G,\Phi)\}/AS)^*/IH^*\cong \fr_S.\eqno{\qed}$$

The deformation retraction of $K_r$ onto the subcomplex spanned by 1-particle irreducible graphs restricts to a
deformation  retraction of $L_S$. This sets up an isomorphism with $IP\overline{\G}_S$, which is the quotient complex of
$P\overline{\G}_S$ by graphs with separating edges (internal or external). 
 As in the $Out(F_r)$ case, this implies there is a chain complex that
computes
$H^*(\Gamma(S);\R)$ which carries a Lie bialgebra structure.  This bialgebra structure vanishes at the level of homology, so that $\oplus_S
H^*(\Gamma(S);\R)$ is the primitive part of a differential graded Hopf algebra.  

\section{The commutative case, revisited}

\subsection{Kontsevich's graph complex}

When $\O$ is the {\it commutative operad}, each vector space $\O[m]$ is
one-dimensional, with a canonical basis element consisting of a rooted tree (not
planar!) with $m$ numbered leaves.   $\Sp[m]$ has a basis consisting of $*_m$.  Thus an
$\O$-graph for the commutative operad is simply an oriented graph, with no additional
structure attached to the vertices. In the graph complex associated to the commutative
operad, the $k$-chains $\G_k$ have as basis oriented graphs with $k$ vertices, and the
boundary operator $\bdry_E$ acts by collapsing edges.  Note that
we may assume that a commutative graph has no edges which are loops, since a graph
automorphism which reverses the orientation on such an edge also reverses the
orientation on the graph, making that graph zero in the chain group.  This chain
complex of oriented graphs was the one studied in \cite{CoV}.  It was noted there that the
subcomplex of one-particle irreducible graphs, on which the bracket and cobracket are
compatible, does not have the same homology as the full complex, unlike the associative
and Lie cases.  We revisit the commutative case  here to give a  geometric
interpretation of the complex and explain the fact that the subcomplex is not quasi-isometric
to the full complex.
We recently discovered a proof that commutative graph homology without separating edges or
separating \emph{vertices} does compute the commutative graph homology.  This is explained in 
\cite{CoGV}.

\subsection{The simplicial closure of Outer space}

Full Outer space $\OX_r$, as opposed to the spine $K_r$, can be described as a union of open simplices,
one for each vertex of 
$K_r$. Recall that a vertex of $K_r$ is a marked graph $(g,X)$ with $k+1$ edges; points of the corresponding 
open  $k$-simplex are given by assigning non-zero lengths to the edges of $X$, subject to
the constraint that the sum of the lengths must equal 1. Gluing instructions are given by adjacency relations 
in $K_r$: an open face is attached when an
edge is allowed to degenerate to length 0, changing the topological type of the graph but not its fundamental group.  

 This union of open simplices is not a simplicial complex, since each
simplex is missing certain faces.  Specifically, the missing faces correspond to 
proper subgraphs which contain cycles; when the edges in such a subgraph are
collapsed to points, the fundamental group of the graph changes and the marked graph is no longer in $\OX_r$.  We can
complete
 $\OX_r$ to a simplicial complex $\CX_r$ by formally adding in these missing
faces.  The action of $Out(F_r)$ extends, and the union
$\s_r$ of the added simplices
 forms an invariant subcomplex of $\CX_r$.  We denote by $C_*(\X_r,\s_r)$ the corresponding relative simplicial chain complex.

There are two actions of $Out(F_r)$ on $\R$, the trivial action and the action via the composition
$Out(F_r)\to Out({\bf Z}^r)=Gl(r,{
\bf Z})\to \R$ given by $\alpha\mapsto \bar\alpha\mapsto \det(\bar\alpha)$. We denote by $\tilde \R$ the reals with the non-trivial action.
We claim that the commutative graph complex can be identified
with
$C_*(\CX_r, \s_r)\otimes_{Out(F_r)} \tilde\R$, so that
Kontsevich's graph homology is isomorphic to the equivariant homology  $H_*^{Out(F_r)}(\CX_r, A; {\tilde\R})$ (see
\cite{Brown}).  

Kontsevich's graph complex $\G$ is a direct sum of complexes $\G^{(r)}$, where 
$\G^{(r)}_k$ has as basis connected, oriented graphs with $k$ vertices and Euler characteristic $1-r$.   

\begin{proposition}\lbl{equivariant}The equivariant homology 
$H_p^{Out(F_r)}(\CX_r,\s_r; \tilde\R) $ is isomorphic to graph homology $H_{p+1-r}(\G^{(r)})$.
\end{proposition}

\begin{proof}  The relative chains $C_p(\CX_r,\s_r)=C_p(\CX_r)/C_p(\s_r)$ have one basis element for each
$p$-simplex in $\CX_r$ which is not entirely contained in $\s_r$.  Such a
$p$-simplex is given by a marked graph $(g,X)$ with $p+1$ edges. The vertices of this
$p$-simplex correspond to the edges $e$ of $X$ (they are obtained by collapsing all edges other than $e$), so that ordering the edges of
$X$ determines an orientation of the corresponding $p$-simplex.  We write $(g,X,e_0,\ldots,e_p)$ to represent the oriented simplex. 

An element $\alpha\in Out(F_r)$ acts on $(g,X,e_0,\ldots,e_p)$ by changing the marking.  Specifically, represent
$\alpha$ by a homotopy  equivalence $f_\alpha\colon R_r\to R_r$ of the standard rose $R_r$; then $(g,X,e_0,\ldots,e_p)\cdot \alpha = (g\circ
f_\alpha,X,e_0,\ldots,e_p)$.  This action is transitive on markings: to send $(g,X,e_0,\ldots,e_p)$ to
$(g',X.e_0,\ldots,e_p)$, choose any homotopy inverse $g^{-1}$ to $g$ and take
$f_\alpha=g^{-1}\circ g'$.  

 The
marking
$g$ induces an isomorphism
$g_*\colon H_1(R_r)\cong
\R^n\to H_1(X)$, thus determining an orientation on $H_1(X)$.
Recall that one of the equivalent definitions of an orientation on a connected graph $X$ is a an orientation of
$\R E(X)\oplus H_1(X)$ (Proposition~\ref{orient}). 
 Thus we have a map 
$$C_p(\CX_r,\s_r)\otimes_{Out(F_r)} \tilde\R \to \G_{p+1-r}$$
 sending $(g, X, e_0,\ldots,e_p)\otimes 1\mapsto (X,or),$
where $or$ is the orientation of $\R E(X)\oplus H_1(X)$ given by the ordering $e_0,\ldots,e_p$ of $E(X)$ together
with the orientation of $H_1(X)$ induced by $g$.

   To see that this map is well-defined, note that the action of $\alpha$ changes orientation iff $\det(\bar\alpha)=-1$.  Since the map is also
surjective, and since the vector spaces in question are of the same (finite) dimension, it is an isomorphism.
 Under this isomorphism, the  boundary map 
$$\bdry\otimes 1\colon C_p(\CX_r,\s_r)\otimes_{Out(F_r)}\tilde\R\to C_{p-1}(\CX_r,\s_r)\otimes_{Out(F_r)}
\tilde\R$$ is identified with the graph homology boundary 
 map
$\bdry_E\colon\G^{(r)}_{p+1-r}\to
\G^{(r)}_{p-r}$. 
\end{proof}

Note that the twisted coefficients account for the $H_1(X)$ part of a graph's orientation. With trivial coefficients,
$H_*^{Out(F_r)}(\CX_r,\s_r; \R) $ is isomorphic to a version of graph homology where the orientation is given by 
ordering the edges only. This homology is considered in a preprint of Bar-Natan and McKay\cite{BN}.

The twisted equivariant homology groups $H_*^{Out(F_r)}(\CX_r,\s_r; \tilde \R)$ can be interpreted in terms
of certain untwisted relative homology groups as follows.  Let $SOut_r$ denote the kernel of the non-trivial action on $\R$, i.e.\
$SOut_r$ is the preimage of $SL(r,\mathbb Z)$ under the natural map $Out(F_r)\to GL(r,\mathbb Z)$.

\begin{theorem}\lbl{longexact}
There is a long exact sequence
\begin{align*}
\cdots\to H_p^{Out(F_r)}&(\CX_r,\s_r; \tilde \R)\to H_p(\CX_r/SOut_r , A_r/SOut_r )\to\\
&\to H_p(\CX_r/{Out(F_r)}, A_r/{Out(F_r)})\to H_{p-1}^{Out(F_r)}(\CX_r,\s_r; \tilde
\R)\to\ldots.\end{align*} 
\end{theorem}

\begin{corollary}\lbl{graphexact}
Let $\O$ be the commutative operad. Then the shifted $\O$-graph complex $P\overline{\G}_{*+1-r}$ fits into a long exact sequence
\begin{align*}
\cdots\to H_p(P&\overline{\G}_{*+1-r})\to H_p(\CX_r/SOut_r , A_r/SOut_r )\to\\
&\to H_p(\CX_r/Out(F_r), A_r/Out(F_r))\to H_{p-1}(P\overline{\G}_{*+1-r})\to\ldots.\end{align*}
\end{corollary}
In other words, graph homology measures the difference between the
homologies of the relative quotients of $(\CX,A)$ by $SOut$ and by
$Out$.

\proof[Proof of Theorem~\ref{longexact}]
We will suppress the subscripts $r$ to streamline the notation.

 Let $\tau\in Out$ be an automorphism with
$\det(\bar\tau)=-1$ (e.g.
$\tau\colon x_1\leftrightarrow x_2$).  

We have  $C_p(\CX/SOut,A/SOut)\cong C_p^+\oplus C_p^-$, where  $C_p^+$
and $C_p^-$ are the $+1$ and $-1$ eigenspaces of the involution induced by
$\tau$.   The short exact sequence
$$1\to C_p^-\to C_p(\CX/SOut,A/SOut)\to C_p^+\to 0$$
gives  rise to a long exact homology sequence
$$\cdots\to H_p(C_*^-)\to H_p((\CX/SOut, A/SOut )\to H_p(C_*^+
)\to H_{p-1}(C_*^-)\to\ldots$$

Now note that  $
C_p(\CX/SOut,A/SOut)\otimes_{\langle\tau\rangle} \R\cong C_p^-,$
implying $$H_p(C_*(\CX,A)\otimes_{Out}\tilde\R) \cong H_p(C_*^-).$$
Furthermore, the map $C_p(\CX/Out, A/Out )\to C_p(\CX/SOut,A/SOut)$ given
by 
$\sigma\mapsto {1\over 2}(\sigma +\tau\sigma)$ is an isomorphism onto
$C_p^+$.  Therefore the above long exact sequence becomes
\begin{align*}
\cdots\to H_p^{Out}(\CX, A;\tilde\R)&\to H_p(\CX/SOut, A/SOut )\to\\
&\to H_p(\CX/Out, A/Out)\to H_{p-1}^{Out}(\CX, A;\tilde\R)\to\ldots.\tag*{\qed}
\end{align*}

\subsection{Homology of the one-particle irreducible subcomplex versus graph homology}

The fact that the subcomplex of one-particle irreducible graphs has the same homology as the full complex for the
associative and Lie operads followed because the $\O$ graph homology in these cases can be identified with
the homology of a space of graphs, which deformation retracts onto the subspace of one-particle irreducible graphs by uniformly shrinking all separating
edges.  If the deformation retraction of Outer space $\OX_r$ extended equivariantly to
$(\CX_r, \s_r)$, then our geometric interpretation would show that graph homology, too, could be computed
using the  subcomplex of
one-particle irreducible graphs (see \cite{Brown}, Proposition 7.3).  However, there is not even a continuous extension of this deformation retraction
to all of $\CX_r$.  The reason is that
$\s_r$ contains simplices which are obtained by  shrinking all {\it non}-separating edges in a graph to points, leaving
a tree.  If one then tries to shrink all of the (separating) edges in the tree to points, one is left with a single vertex,
which  does not correspond to a point in $\CX_r$.  

Let $\OY_r$ denote the subspace of $\OX_r$ corresponding to 1-particle irreducible graphs. Then the homology of  the subcomplex of the
graph complex spanned by 1-particle irreducible graphs can be identified with the equivariant homology $H_*^{Out(F_r)}(\CY_r, A_r\cap
\CY_r);\tilde\R)$, as in the proof of Proposition~\ref{equivariant}, measuring the difference between the relative quotients of $(\CY_r, A\cap \CY_r)$ modulo
$Out(F_r)$ and modulo
$SOut(F_r)$.   Although $H_*^{Out(F_r)}(\CY_r, A_r\cap
\CY_r;\tilde\R)\cong H_*^{Out(F_r)}(\CX_r, A_r;\tilde \R)$ for small values of $r$, this is accidental, as the geometry suggests; for larger values of $r$
the homologies differ.  We illustrate the difference between $\CX$ and $\CY$ for $r=2$, where one can easily draw a picture, below.
 
For $r=2$, $\OY_2$ can be identified with the familiar picture of the
hyperbolic plane, tiled by ideal triangles (Figure~\ref{Y}).
\begin{figure}[ht!]
\begin{center}
\includegraphics[height=5cm]{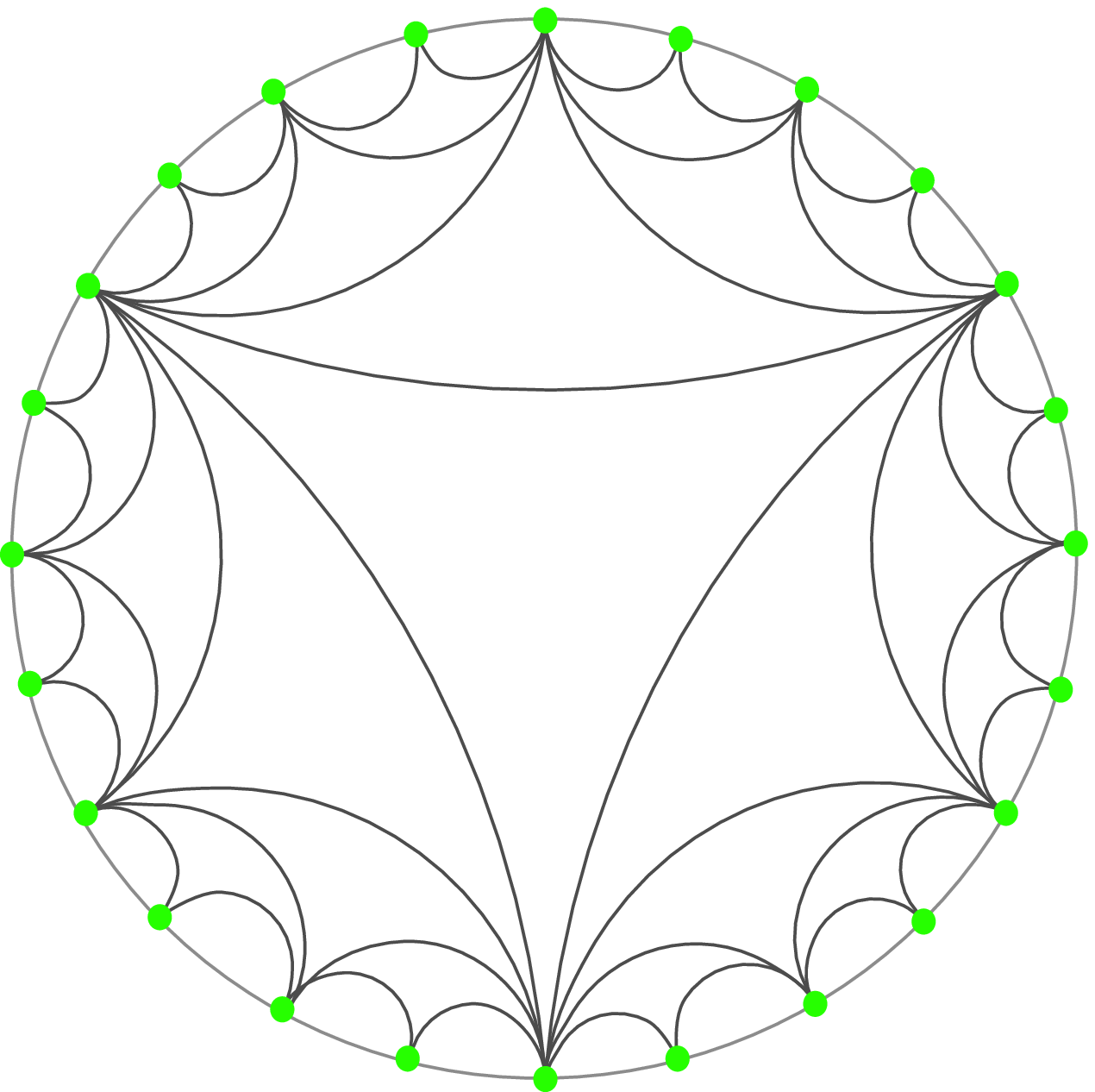} 
\end{center}
\caption{The space $\CY_2$}\lbl{Y}
\end{figure}
Full Outer space $\OX_2$ has one additional
open triangle (a ``fin") adjacent to each edge in the ideal triangulation. (Figure~\ref{X}).
\begin{figure}[ht!]
\begin{center}
\includegraphics[height=3cm]{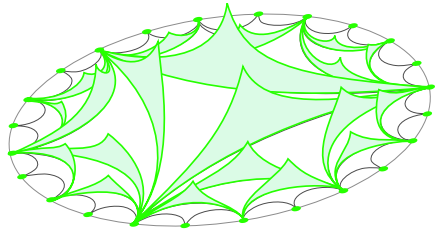} 
\end{center}
\caption{Outer space in rank 2}\lbl{X}
\end{figure}
The simplicial closure $\CY_2$ is obtained by adding the vertices of the ideal
triangles; the closure of $\CX_2$ contains also the additional vertex and two missing edges of each fin.

We conclude by using this example to illustrate Corollary~\ref{graphexact}. The quotient of $\CX_2$ by the action of $SOut_2$ is a  
``pillow" (perhaps more familiar as the quotient of the hyperbolic plane by the action of $SL(2,\mathbb Z)\cong SOut_2$) with half of a fin
attached (Figure~\ref{pillow});
\begin{figure}[ht!]
\begin{center}
\includegraphics [height=3.5cm]{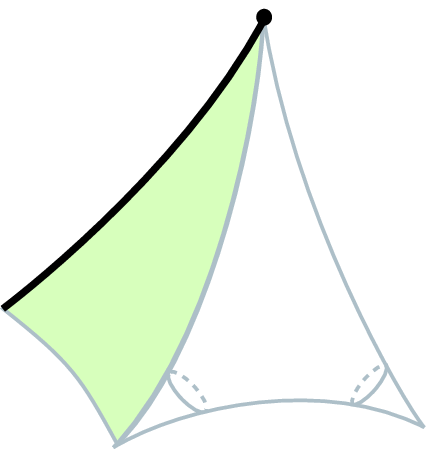}
\end{center}
\caption{Quotient of $\CX_2$ modulo $SOut_2$}\lbl{pillow}
\end{figure}
 the image of $A_2$ is a (closed)
edge of this fin.  Therefore the relative homology
$H_*(\CX_2/SOut_2,A_2/SOut_2)$ has one generator in dimension 2 and is zero
elsewhere.  In the quotient by the full outer automorphism group $Out(F_2)$, the pillow is flattened to a triangle; again the image of $A_2$
is an edge of the attached fin, so that the relative homology vanishes.  Therefore the only non-zero terms
in the long exact sequence are
$$\cdots 0\to H_2^{Out(F_2)}(\CX_2,A_2)\to \R\to 0 \to \cdots.$$
Proposition~\ref{equivariant} then tells us that $H_2(\G^{(2)})\cong \R$. In terms of oriented graphs, the graph homology generator 
is the theta graph, with two vertices connected by three edges.

\Addresses\recd
\end{document}